\documentclass[11pt]{article}

\usepackage[a4paper, total={16cm,25cm}]{geometry}

\usepackage[style=alphabetic,
  backend=biber,
  maxnames=3,
  maxbibnames=99,
  backref,
]{biblatex}
\usepackage{enumitem}
\setenumerate{label=\textup{(\roman*)}}

\usepackage{mlmodern}

\usepackage{tabularx, booktabs, diagbox}

\usepackage{amsmath, mathtools}
\usepackage{amssymb}
\usepackage{amsthm, thmtools}
\usepackage{tikz-cd}
\usepackage{braket}
\usepackage{bm}

\usepackage[english=american]{csquotes}

\usepackage{todonotes}

\usepackage[colorlinks]{hyperref}
\usepackage[nameinlink]{cleveref}
\makeatletter
\AddToHook{cmd/appendix/before}{\def\cref@section@alias{appendix}}
\makeatother

\declaretheorem[numberwithin=section]{theorem}
\declaretheorem[numberlike=theorem]{lemma, corollary, proposition}
\declaretheorem[numberlike=theorem,style=definition]{definition}
\declaretheorem[numberlike=theorem,style=remark]{example}
\declaretheorem[numbered=no, name=Theorem]{tool}
\declaretheorem[numberlike=theorem,style=remark]{remark}
\declaretheoremstyle[
notefont=\bfseries, notebraces={(}{)}
]{conjecture}

\newcommand{\RR}{\mathbb{R}}

\newcommand{\CC}{\mathbb{C}}

\newcommand{\PP}{\mathbb{P}}

\DeclareMathOperator{\rank}{rank}

\usepackage{xspace}

\usepackage{adjustbox}
\newcommand{\Rset}{\ensuremath{\mathbb{R}}}
\newcommand{\Cset}{\ensuremath{\mathbb{C}}}

\newcommand{\Nset}{\ensuremath{\mathbb{N}}}
\newcommand{\Pset}{\ensuremath{\mathbb{P}}}
\DeclareMathOperator{\image}{image}
\DeclareMathOperator{\codim}{codim}
\newcommand{\h}{{\mathrm h}}
\newcommand{\nullity}[1]{\ensuremath{\operatorname{nullity} #1}}
\newcommand*{\tran}{^{\mkern-1.5mu\mathsf{T}}}
\DeclarePairedDelimiter{\norm}{\lVert}{\rVert}
\newcommand{\vect}[1]{\ensuremath{\boldsymbol{#1}}}
\newcommand{\matx}[1]{\ensuremath{\boldsymbol{#1}}}
\newcommand{\Hardy}{\ensuremath{\mathcal{H}_2}}
\newcommand{\stabModels}[1]{\ensuremath{\mathcal{M}_{#1}}}
\newcommand{\stabDiscModels}[1]{\ensuremath{\mathcal{M}_{#1,\mathrm{d}}}}
\newcommand*\diff{\mathop{}\!\mathrm{d}}
\newcommand{\bmat}[1]{\ensuremath{
    \begin{bmatrix}
        #1
    \end{bmatrix}
}}
\newcommand{\pmat}[1]{\ensuremath{
    \begin{pmatrix}
        #1
    \end{pmatrix}
}}
\newcommand{\Toeplitz}[2]{\ensuremath{\matx{T}_{#2}(#1)}}

\newcommand{\qedex}{\hfill $\triangle$}

\title{Porteous Meets Walsh: Counting Critical Points for \texorpdfstring{$\Hardy{}$}{H₂}-Norm Model Order Reduction Using Algebraic Geometry Behind Multiparameter Eigenvalue Problems}
\date{}

\usepackage{authblk}
\author[1]{Sibren Lagauw}
\author[2]{Leonie Kayser}
\author[1,3]{Bart De Moor}
\affil[1]{STADIUS-ESAT, KU Leuven, Leuven.}
\affil[2]{MPI for Mathematics in the Sciences, Leipzig.} 
\affil[3]{Fellow IEEE, IFAC, SIAM.}

\usepackage{xcolor}

\begin{document}
\maketitle

\begin{abstract}
    Several polynomial systems arising in systems and control have been shown to exhibit a partially linear structure, which allows them to be reformulated as a multiparameter eigenvalue problem (MEP).
    In this work, we show how such MEP formulations can be leveraged to obtain refined information about the solution set of the polynomial system.
    We demonstrate, in the setting of single-input single-output $\Hardy$-norm model order reduction, how the Thom--Porteous formula can be used to derive an exact count of the complex solutions of the critical equations derived from Walsh's theorem.
    We prove, under suitable genericity assumptions, that their set of complex solutions is finite and that all solutions occur without multiplicity.
    Moreover, we derive an explicit formula for the total number of complex solutions as a function of the model orders.
    Throughout, we provide a self-contained and application-oriented exposition of the algebro-geometric tools required for our analysis, with the aim of making them accessible to the broader engineering community.
\end{abstract}

\section{Introduction}
\label{sec:intro}
\addcontentsline{toc}{section}{\nameref{sec:intro}}

Systems of polynomial equations arise frequently in engineering applications; still, solving them comprehensively---numerically or symbolically---is computationally challenging.
For such systems, one typically has access to unconditional but only coarse upper bounds on the number of affine solutions, such as those provided by Bézout’s theorem.
Nevertheless, having a precise a priori count of the solutions is valuable both for the design and use of computational methods. 
It provides a concrete proxy for the expected computational cost and can be used to validate computations by confirming that the full solution set has been recovered.

A \emph{multiparameter eigenvalue problem} (MEP) generalizes the standard eigenvalue problem $(\matx{A} - \lambda \matx{I})\vect{v} = \vect{0}$ by allowing multiple spectral parameters and by dropping the restriction of linearity in the parameters~\cite{Atkinson1972,Volkmer1988,DeCock2021}.
For rectangular MEPs~\cite{Vermeersch2023,Hochstenbach2024}, one looks for all $\vect{\lambda} = (\lambda_1, \dots, \lambda_k) \in \mathbb{C}^k$ and nonzero $\vect{v} \in \mathbb{C}^m$ such that 
\begin{equation} \label{eq:mep-def}
    \matx{M}(\vect{\lambda}) \vect{v} = \vect{0},
\end{equation}
where $\matx{M}(\vect{\lambda})$ is a tall, $\ell \times m$ matrix whose entries are polynomials in $\vect{\lambda}$.  
The solutions $\vect{\lambda}$ are called \emph{eigentuples} \cite{DeCock2021} and $\vect{v}$ the corresponding eigenvector. 
For $k = 1$ and $\matx{M}$ linear in $\lambda$, we recover the (generalized) eigenvalue problem.

Several systems of polynomial equations arising in applications from systems and control have been shown to exhibit a partially linear structure, such that they can be reformulated as a rectangular MEP, see for example \cite{Batselier2012,Vermeersch2019,Moor2020,Agudelo2021,Lagauw2025}.
In these cases, instead of solving the system directly using general-purpose polynomial system solvers (e.g., \cite{Breiding2018}), one can alternatively obtain the solutions by computing the eigentuples of the corresponding MEP using specialized solvers (e.g., \cite{Vermeersch2023,Hochstenbach2024,vermeersch2025a}).

In this paper, we argue that the reformulation of the system of polynomial equations as an MEP can also yield more theoretical structural insights into the problem at hand.
In particular, two central results in algebraic geometry, when applied in the context of these MEP reformulations, can be used to obtain finer information about the solution set of the polynomial system.
When a polynomial system is expected to have finitely many solutions, the \emph{fiber dimension theorem} is a useful tool to verify this. 
Consequently, where finiteness of the solutions can be established, the \emph{Thom--Porteous formula}, applied to the determinantal variety associated with the MEP reformulation, can often be used to derive an exact count of the affine complex solutions.

To ground our discussion, we illustrate the use of these tools in the context of a problem from systems and control: the \Hardy{}-norm model reduction problem for single-input single-output (SISO) linear time-invariant (LTI) models, in which a given model is approximated by one of lower order. 
It corresponds to a non-convex optimization problem and may therefore admit multiple local optima, which can be characterized as (a subset of the) real-valued common roots of a system of polynomial equations~\cite{DeMoor1993,Regalia1995,Hanzon2007,Lagauw2023}.
These equations, commonly referred to as the \emph{critical equations}, follow from \emph{Walsh's Theorem}~\cite{Walsh1960}, which, together with the Thom--Porteous formula, motivates the title of this paper.

We study the complex solution set of the critical equations as a whole and refer to its elements as the \emph{critical points}.
When the given higher-order model is sufficiently generic, the set of critical points is known to be finite~\cite[Theorem~3]{Baratchart1986}; in fact, simple poles suffice for this to hold~\cite[Theorem~3.2]{Hanzon2012}.
This observation has motivated several approaches that use algebraic and numerical linear algebra methods to enumerate all such critical points; see, e.g., \cite{Hanzon2007,Peeters2012,Alsubaie2019,Agudelo2021,Lagauw2023}.
When complete, such enumerations allow one to certify global optimality.  
However, they are computationally feasible only for relatively small problems and are therefore best regarded as theoretical tools, mostly used for benchmarking purposes. 
Practical solution approaches, such as the iterative rational Krylov algorithm (IRKA), scale to large problems by iteratively enforcing the first-order optimality conditions, but converge at best to a single local optimum, without any guarantee of global optimality; see, e.g., \cite{Gugercin2008} and references therein. 
Despite the extensive literature on $\mathcal{H}_2$-norm model reduction, an exact count of the critical points for arbitrary problem sizes is, to the best of our knowledge, not available. For the co-order one case, i.e., when the order is reduced by one, \cite{Hanzon2007} provides an upper bound.

Recent work showed that the critical equations allow for an MEP reformulation~\cite{Lagauw2023}. 
By applying the algebro-geometric tools outlined above, we arrive at the following result:

\begin{tool}[\Cref{thm:numberSolutions}]
For generic higher-order models of McMillan degree $N$, the number of critical points of the model order reduction problem to degree $n$, $1\leq n<N$, is
\[
\sum_{j=0}^n \binom{N-n-1+j}{j} 2^j.
\]
\end{tool}
The main contributions of our present work are as follows:
\begin{enumerate}[label=(\arabic*)]
    \item We give a new proof that, for higher-order models with simple poles, the set of critical points is finite (\Cref{cor:finiteness_by_simple_roots}). 
    Moreover, we show for generic higher-order models that the solutions of the critical equations occur without multiplicity, a property favorable for numerical computation (\Cref{cor:finiteness}).

    \item Assuming the higher-order model has only simple poles, we derive an explicit formula for the number of critical points, that is, the total number of complex solutions to the critical equations, as a function of the full and reduced orders (\Cref{thm:numberSolutions}). 
    For higher-order models with one or more multiple poles, we illustrate how the aforementioned count can be used to derive an upper bound on the number of critical points (\Cref{sec:multiplepoles}). 
    Although we consider the continuous-time case throughout, our analysis easily extends to a discrete-time setting (\Cref{sec:discretetime}).

    \item We provide a self-contained, application-oriented exposition of the key algebro-geometric tools, making them accessible to a broader engineering audience (\Cref{sec:finitAndDeg} and \Cref{app:AG}). These techniques are not specific to the model reduction problem, but are equally applicable to other engineering problems where such MEP formulations arise.
\end{enumerate}

We would like to note that our results connect to a broader line of work in algebraic geometry.
Counting of critical points is related to the \emph{Euclidean distance (ED) degree} of an algebraic variety, defined as the number of complex critical points of the squared Euclidean distance from a generic data point to the smooth part of the variety~\cite{Draisma2016}. 
In discrete time, the squared $\mathcal{H}_2$-norm of a SISO LTI system is the sum of its squared impulse response coefficients. 
The discrete-time model reduction problem can be seen as an infinite-dimensional instance of structured low-rank approximation of \emph{Hankel matrices}, whose ED degrees are studied in~\cite{Ottaviani2014,Panah2025MiddleCatalecticant} and \cite{Kayser2026a}.
The latter employs similar algebro-geometric tools as used here, and can be seen as a companion to the present paper.

\paragraph{Outline} 
The remainder of the paper is organized as follows.
In \Cref{sec:CritEquations}, we give a formal statement of the model reduction problem and derive the critical equations of the optimization problem.
The heart of the paper is \Cref{sec:finitAndDeg}, where we lay out a recipe for establishing finiteness and deriving a count of the solution set of MEPs.  
We implement this recipe for the model reduction problem and derive the corresponding solution count.
Several further considerations are given in \Cref{sec:otherConsid}: the special case where the given higher-order model has multiple poles and the discrete-time setting. 
We illustrate all theoretical results by means of several examples in \Cref{sec:Numerical} and conclude the paper in \Cref{sec:Conclusions}.
The complementary \Cref{app:AG} gives an elementary introduction to the essential algebraic geometry concepts relevant to our work.

\paragraph{Notation}
Throughout, scalars are denoted by lowercase letters (e.g., $s\in\Cset$), vectors by boldface lowercase letters (e.g., $\vect{a}\in\Cset^{k}$), and matrices by boldface uppercase letters (e.g., $\matx{M} \in \Cset^{k \times k}$).
Let $\CC[s]_{\leq k}$ denote the $(k+1)$-dimensional vector space of complex polynomials in $s$ of degree at most $k$. 
For a polynomial $f=\sum_{i=0}^{k}f_{i}s^{i} \in \CC[s]_{\leq k}$, we write its coefficient vector in boldface $\vect{f}=(f_{0},\dots,f_{k})\tran \in \Cset^{k+1}$, and conversely, any vector $\vect{f} \in \Cset^{k+1}$ determines a polynomial $f \in \CC[s]_{\leq k}$.
Let $\Toeplitz{f}{\ell} \in \Cset^{\ell \times (\ell + k)}$ denote the banded Toeplitz matrix generated by $\ell$ forward-shifts of the coefficient vector $\vect{f} \in \Cset^{k+1}$:
\begin{equation*}
\Toeplitz{f}{\ell} \coloneq \bmat{
f_0  & \cdots & f_{k-1} & f_k & 0 & \cdots & 0 \\
0 & f_0  & \cdots & f_{k-1} & f_k & \cdots & 0 \\
\vdots  & \ddots & \ddots  & \ddots  & \ddots  & \ddots & \vdots \\
0  & \cdots & 0 & f_0  & \cdots & f_{k-1} & f_k
}.
\end{equation*} 
Toeplitz matrices constructed in this way allow us to represent polynomial multiplication in terms of coefficient vectors. For instance, if $h \in \CC[s]$, then the product $hf$ has coefficient vector $\Toeplitz{h}{k+1}\tran \vect{f}$. 
For $f \in \mathbb{C}[s]_{\le k}$, write $f^\h \coloneq  t^k f(s/t) \in \CC[s,t]_k$ for its homogenization to degree $k$.
Algebraic varieties are denoted by uppercase letters, with calligraphic font occasionally used for emphasis (e.g., $\mathcal{W}$).
Unless stated otherwise, all open and closed subsets mentioned in this paper are taken with respect to the Zariski topology.

\section{\texorpdfstring{$\Hardy{}$}{H₂}-norm model order reduction}\label{sec:CritEquations}

Throughout, we work with transfer functions of real, continuous-time, strictly causal, stable, single-input single-output (SISO) linear time-invariant (LTI) systems of finite order (but see \Cref{sec:otherConsid} for the discrete-time setting).
For a fixed integer $N \geq 1$, let $\stabModels{N}$ denote the set of all such transfer functions of \emph{McMillan degree} exactly $N$.
Any model $H \in \stabModels{N}$ can be parametrized as  
\[
H(s) = \frac{d}{c} = \frac{d_{N-1} s^{N-1} + \dots + d_{1} s + d_0}{c_Ns^N + c_{N-1} s^{N-1} + \dots + c_{1} s + c_0},
\]
where $c \in \Rset[s]_{\leq N}$ and $d \in \Rset[s]_{\leq N-1}$ are real polynomials with $c_N \neq 0$ and $\gcd(c,d)=1$.
\emph{(Hurwitz-)Stability} here means that all poles of $H$, the roots of $c$, lie in the open left half-plane.
\emph{Strict causality} means $\deg c > \deg d$, so that the impulse response vanishes for negative time and is square-integrable for $t \ge 0$.
Consequently, elements of $\stabModels{N}$ are real rational functions in the \emph{Hardy space} $\Hardy$ of functions analytic in the open right half-plane and square-integrable along the imaginary axis.

Given a target reduced order $1 \leq n < N$, we seek an approximant $\widehat{H} \in \stabModels{n}$, parametrized as
\[
\widehat{H}(s) = \frac{b}{a} = \frac{b_{n-1} s^{n-1} + \dots + b_{1} s + b_0}{a_n s^n + a_{n-1} s^{n-1} + \dots + a_{1} s + a_0},
\]
where $a\in\Rset[s]_{\leq n}$ and $b \in \Rset[s]_{\leq n-1}$ are real polynomials with $a_n \neq 0$, $\gcd(a,b) = 1$, and $a$ Hurwitz stable. 
The goal of the model order reduction problem under consideration is to minimize the $\Hardy$-norm of the \emph{approximation error} $H - \widehat{H}$.
Formally, the problem corresponds to
\begin{align}
\operatorname*{minimize}_{\widehat{H} \in \stabModels{n}} \, J^2, 
\qquad\qquad
J^2 = \norm{H - \widehat{H}}_{\Hardy{}}^2 &= \frac{1}{2\pi} \int_{-\infty}^{\infty} |H(\mathrm{j}\omega) - \widehat{H}(\mathrm{j}\omega)|^2 \diff \omega \label{eq:modelReductionOptProblem} \\
&= \int_0^{\infty} |h(t) - \widehat{h}(t)|^2 \diff t,\nonumber
\end{align}
where $h(t)$ and $\widehat{h}(t)$ denote the impulse responses of $H$ and $\widehat{H}$, respectively, obtained via the inverse Laplace transform.

    We restrict the search to reduced models of McMillan degree exactly $n$, rather than at most $n$. 
    This is because the set $\mathcal{M}_n$ is open and dense in $\mathcal{M}_{\le n} \coloneq \mathcal{M}_0 \cup \mathcal{M}_1 \cup \dots \cup \mathcal{M}_n$, where $\mathcal{M}_0 \coloneq \{0\}$, so that every model of lower McMillan degree is a limit of models of degree $n$. 
    Moreover, every local minimizer of \eqref{eq:modelReductionOptProblem}, when it were to be considered over $\mathcal{M}_{\le n}$, is guaranteed to have McMillan degree exactly $n$~\cite[Theorem~5.3]{Regalia1995}.\footnote{The result is stated in~\cite{Regalia1995} for discrete time and is proven by showing that when a reduced model of lower McMillan degree is perturbed by adding a first-order term, then there always exists a perturbation for which the approximation error of the perturbed model decreases. A completely similar argument can be made in the present, continuous-time setting.} 
    The set $\mathcal{M}_n$ is a smooth manifold of dimension $2n$: indeed, normalize $a$ to be monic so that $a_n=1$ and the remaining coefficients of $a$ and $b$ range over an open subset of $\mathbb{R}^{2n}$, obtained by removing those points where $a$ and $b$ share one or more common factors. 
    We call $\widehat{H} \in \mathcal{M}_n$ a \emph{stationary point} of~\eqref{eq:modelReductionOptProblem} if the gradient of $J^2$ with respect to these $2n$ coefficients vanishes.


The first-order optimality conditions of \eqref{eq:modelReductionOptProblem} can be derived in two complementary ways: analytically, by partial differentiation of the $\Hardy{}$-norm error~\cite{Aigrain1949,Meier1967}, and geometrically, by exploiting orthogonality in the Hardy space $\Hardy$ (see, e.g., \cite{Regalia1995,Gugercin2008}).
By doing so, one obtains a set of Hermite interpolation constraints at the reflected poles of the reduced model $\widehat{H}$.
This foundational result, originally due to Aigrain and Williams~\cite{Aigrain1949}, was established for reduced models with distinct, real poles and later extended by Meier and Luenberger~\cite{Meier1967} to cover complex and higher-order poles.

\begin{theorem}[Meier \& Luenberger~\cite{Meier1967}]\label{thm:interpolation}
    Consider, for a given $H \in \stabModels{N}$, the model reduction problem as formulated in \eqref{eq:modelReductionOptProblem}. 
    A necessary condition for a reduced model $\widehat{H} \in \stabModels{n}$ to be a local minimizer is that, for each pole $\hat p_i \in \Cset$ of $\widehat{H}$, the following interpolation conditions hold:
    \begin{equation}
        H^{(k)}(- \hat p_i) = \widehat{H}^{(k)}(- \hat p_i), \quad k = 0, \dots, m_i, \label{eq:interpolationConditions_Cont}
    \end{equation} where $m_i$ is the pole order of $\hat p_i$ and $F^{(k)}(z)$ denotes the $k$th derivative of $F$, evaluated at $z$.
\end{theorem}

Model reduction for linear systems is closely related to the mathematical discipline of function approximation.
From this viewpoint it is unsurprising that similar optimality conditions were already known in the field of classical analysis in the 1930s, e.g., in work by Walsh~\cite{Walsh1960}. 
Adapted to our context of model reduction, the result reads as follows.

\begin{theorem}[Walsh's theorem]\label{H2MOR_thm:walsh}
Consider a higher-order model $H \in \mathcal{M}_{N}$. 
An approximant $\widehat{H} \in \mathcal{M}_{n}$ is a stationary point of the model reduction problem~\eqref{eq:modelReductionOptProblem} if and only if, as functions,
\begin{equation*}
    H - \widehat{H} = \frac{d}{c} - \frac{b}{a} = \left( \frac{\tilde{a}}{a} \right)^2 R, \label{H2MOR_eq:Walsh}
\end{equation*} with $\tilde a \coloneq \sum_{i=0}^n a_i(-s)^i \in \Rset[s]_{\leq n}$, which has the mirror images $-\hat p_i$ of the poles $\hat p_i$ as its roots with the same multiplicities, and $R$ a real rational function in $\mathcal{H}_2$.
\end{theorem} 

\Cref{H2MOR_thm:walsh} is the continuous-time counterpart of~\cite[Theorem 3.14]{Regalia1995}.
It can be derived using a projection argument in the Hardy space $\mathcal{H}_2$~\cite[Section~II]{Hanzon2007},~\cite[Section~8.3]{Bleylevens10}. 
Similar results have been obtained using state-space realizations~\cite[Theorem~4.5]{VanDoren2010}.

\begin{proof}[Proof sketch]
Since $\mathcal{M}_n$ is a smooth submanifold of $\mathcal{H}_2$, the model $\widehat{H}$ is a stationary point if and
only if the approximation error $H - \widehat{H}$ is orthogonal to the tangent space of $\mathcal{M}_n$ at $\widehat{H}$.

Normalize $(a,b)$ so that $a$ is monic. 
Differentiating $\widehat{H} = b/a$ along a perturbation $(a + \epsilon\alpha, b + \epsilon\beta)$, with $\alpha, \beta \in \mathbb{R}[s]_{\le n-1}$, gives the tangent vector
\[
  \frac{\mathrm{d}}{\mathrm{d}\epsilon}\, \frac{b + \epsilon\beta}{a + \epsilon\alpha}
  \bigg|_{\epsilon = 0} = \frac{a\beta - b\alpha}{a^2}.
\]
Since $\gcd(a,b) = 1$, the linear map $(\alpha, \beta) \mapsto a\beta - b\alpha$ is injective, and hence, by a dimension count, onto $\mathbb{R}[s]_{\le 2n-1}$. 
The tangent space of $\mathcal{M}_n$ at $\widehat{H}$ therefore consists of the rational functions $f/a^2$ with $f \in \mathbb{R}[s]_{\le 2n-1}$, regardless of the multiplicities of the poles of $\widehat{H}$. 
Over $\mathbb{C}$, it is spanned by the functions $1/(s - \hat p_i)^k$, $k = 1, \dots, 2m_i$, where $m_i$ is the order of the pole $\hat{p}_i$.

The function $(\tilde a/a)^2$ is analytic in the open right half-plane and has modulus one on the imaginary axis, and by the Beurling--Lax theorem~\cite{Lax1959}, see also~\cite[Theorem 3.1.5]{Partington2004}, the orthogonal complement in $\mathcal{H}_2$ of functions $f/a^2$ is the subspace $(\tilde a/a)^2 \mathcal{H}_2$.
Hence the approximation error $H - \widehat{H}$ is orthogonal to the tangent space if and only if
$H - \widehat{H} = (\tilde a/a)^2 R$ for some $R \in \mathcal{H}_2$. Clearly, $R$ must be a real rational function as well.
\end{proof}

\Cref{H2MOR_thm:walsh} implies that the numerator of $H - \widehat{H} = (ad - bc)/(ac)$
satisfies 
\begin{equation}
    ad - bc = \tilde a^2 R c/a. \label{eq:walshnodenom}
\end{equation}
The zeros of $\tilde a^2$ lie in the open right half-plane, whereas the poles of $R c/a$, which are poles of $R$ or zeros of $a$, lie in the open left half-plane by the stability assumption, so that none of the zeros of $\tilde a^2$ can be canceled.
Hence, $\widehat{H} \in \mathcal{M}_n$ is a stationary point if and only if $\tilde a^2$ divides $ad - bc$, or, equivalently, if and only if the interpolation conditions~\eqref{eq:interpolationConditions_Cont} hold for $k = 0, \dots, 2m_i - 1$. See \Cref{rem:multiplicity} at the end of this section on the subtle differences between \Cref{thm:interpolation,H2MOR_thm:walsh}.

The left-hand side of \eqref{eq:walshnodenom} is polynomial of degree $N+n-1$, implying that $\widehat{H} \in \mathcal{M}_n$ is a stationary point if and only if there exists a polynomial $g \in \Rset[s]_{\leq N-n-1}$ such that
\begin{equation}
  ad - bc = \tilde a^2 g. \label{eq:lofs}
\end{equation}
Since \eqref{eq:lofs} must hold identically as a polynomial identity in $s \in \mathbb{C}$, it is equivalent to the vanishing of each coefficient in the expansion of the polynomial $ad-bc-\tilde{a}^2g$.
Enforcing these coefficients to vanish gives a system of $N+n$ polynomial equations in the unknown coefficients of the polynomials $a,b$, and $g$. 
It is equivalently given by the following matrix expression:
\begin{equation}
 \Toeplitz{a}{N}\tran \vect{d} - \Toeplitz{c}{n}\tran \vect{b} - \Toeplitz{\tilde{a}}{N}\tran \Toeplitz{\tilde{a}}{N-n}\tran \vect{g} = \vect{0}. \label{eq:criticalEqsMatrix}
\end{equation}

We refer to this system of $N+n$ polynomial equations in the $N+n+1$ unknown coefficients $a$, $b$, and $g$ as the \emph{critical equations}.
Since the system is invariant under scaling $(a,b,g) \mapsto (\lambda a, \lambda b, \lambda^{-1} g)$ for $\lambda \in \Cset \setminus \{0\}$, one can expect finitely many solutions up to normalizing (e.g., by $a_n=1$).
We call the complex solutions with $a \neq 0$, considered up to scaling, the \emph{critical points}. 
By~\Cref{H2MOR_thm:walsh}, real critical points with $\deg a = n$, $\gcd(a,b) = 1$, and $a$ Hurwitz stable are exactly the stationary points of~\eqref{eq:modelReductionOptProblem}; and in particular, they include all local minimizers.

Note that the coefficients of $b$ and $g$ appear only linearly in \eqref{eq:criticalEqsMatrix}. 
It is precisely this partially linear structure that is exploited in the corresponding MEP formulation (cf.~\Cref{sec:MEP}).

\begin{remark}[Multiplicity structure]\label{rem:multiplicity}
    For approximants with simple poles, Walsh's theorem enforces exactly the interpolatory conditions of \Cref{thm:interpolation}. 
    However, for a pole $\hat p_i$ of multiplicity $m_i \geq 2$, it imposes $m_i - 1$ additional conditions.
    The difference in both formulations stems from the parametrization considered. 
    Meier and Luenberger~\cite{Meier1967} expand $\widehat{H}$ in partial fractions with a prescribed multiplicity structure: each distinct pole of multiplicity $m_i$ is described by a single pole location and $m_i$ residues. 
    This parametrization only covers the models with that particular multiplicity structure, and therefore covers only a subset of $\mathcal{M}_n$ of dimension $\sum_i (m_i + 1)$, where the sum runs over the distinct poles, which is smaller than $2n$ as soon as one of the poles is multiple. 
    \Cref{thm:interpolation} describes the stationarity conditions on this submanifold, i.e., the first-order optimality conditions for the model reduction problem constrained to approximants with the prescribed multiplicity structure. 
    The $m_i - 1$ additional conditions in Walsh's theorem correspond to perturbations that split the multiple
    pole, and thereby leave this submanifold.
    Consequently, every stationary point of~\eqref{eq:modelReductionOptProblem} satisfies the conditions of \Cref{thm:interpolation}, but a reduced model with multiple poles may satisfy them without being a stationary point of~\eqref{eq:modelReductionOptProblem}. 
\end{remark}

\section{Finiteness and degree of the critical points}\label{sec:finitAndDeg}

In our analysis of the critical points, we deal with several products of vector spaces. 
To ease notation, we use the following shorthands:
\[
A \coloneq \CC[s]_{\leq n}, \quad B\coloneq \CC[s]_{\leq n-1}, \quad C \coloneq \CC[s]_{\leq N}, \quad D\coloneq \CC[s]_{\leq N-1}, \quad G \coloneq \CC[s]_{\leq N{-}n{-}1}.
\]
Additionally, $A^\circ \coloneq A \setminus 0$, and we use, for example, $\mathit{A^\circ\!BG}$ to denote the product $(A\setminus 0) \times B \times G$.
Consider the incidence variety consisting of tuples $(a,b,c,d,g)$ satisfying the critical equations.

\begin{definition}
The \emph{Walsh variety} is $\mathcal{W} \coloneq \set{ (a,b,c,d,g) | a d - b c = \tilde{a}^2 g } \subseteq \mathit{A^\circ\! BCDG}$.
\end{definition}

In this section, we use the Walsh variety to derive an exact count of the critical points.
Our approach relies on tools from algebraic geometry; to keep the exposition self-contained and accessible to a broader engineering audience, the required background is presented in \Cref{app:AG}.
The overall strategy can be summarized as follows:

\begin{itemize}[wide]
\item We first analyze the incidence variety $\mathcal{W}$, whose variables naturally split into two groups: the \emph{input variables} ($c,d$), specifying a problem instance, and the \emph{unknown variables} ($a,b,g$), representing candidate solutions. 
While the fibers of the variety over fixed input data are difficult to analyze directly, the variety admits a simpler description after fixing a different subset of variables (see, e.g., \Cref{lem:divisibilityConstraint} and the proof of \Cref{thm:V_sys_finite}). 
This then allows us to derive key geometric properties of the incidence variety, such as its dimension, irreducibility, and reducedness (\Cref{thm:V_sys_finite}).
\item Once we show that the dimension of the incidence variety coincides with the dimension of the input space, the fiber dimension theorem (\Cref{thm:fiber_dimension}) implies that for generic input data the corresponding fiber consists of (at most) finitely many isolated points (\Cref{cor:finiteness}).
Moreover, if the incidence variety was shown to be defined in a reduced way, then the fiber dimension theorem implies that these generic fibers are reduced as well.
\item To obtain an explicit solution count, we then exploit the partially linear structure in the system of equations by passing to its MEP formulation. 
This yields a determinantal description in which the variables that appear only linearly ($b,g$) are eliminated, at the cost of introducing potential extraneous solutions.
These spurious points form what we call the \emph{bad locus}. 
We show that for generic input data the bad locus is empty (\Cref{lem:MORbmatrix}), so that the determinantal variety is in fact isomorphic to the generic fibers of the incidence variety over the input data (\Cref{cor:equivalenceMEP}).
Moreover, by \Cref{lem:reducedMinors}, reducedness of these fibers implies reducedness of the determinantal variety.
Finally, we apply the Thom--Porteous formula (\Cref{thm:TP}) to compute the degree of the determinantal variety, obtaining the desired number of solutions (\Cref{thm:numberSolutions}).
\end{itemize}

\begin{remark}
    We have deliberately formulated the above strategy in a general way since it applies beyond the present setting of model reduction, and can be used for other systems of equations and MEPs arising in related applications.
    The only problem-specific ingredients are the analysis of the incidence variety (cf.~\Cref{thm:V_sys_finite}) and of the bad locus (cf.~\Cref{lem:MORbmatrix}). 
    See, for instance, our related work~\cite{Kayser2026a} on Euclidean distance degrees of secant varieties of the rational normal curve. 
    In contrast to the present setting, reducedness, irreducibility, and the dimension of the incidence variety are relatively easily established in that context, whereas the corresponding bad locus is substantially more intricate to characterize than the one encountered here.
\end{remark}

\subsection{Fiber dimension theorem and Thom--Porteous}\label{sec:AGtools}

Before implementing the above strategy on the incidence variety $\mathcal{W}$, we recall from \Cref{app:AG} two classical, yet powerful, algebro-geometric tools that play a central role in our analysis: the fiber dimension theorem and a special case of the Thom--Porteous formula.

The \emph{fiber dimension theorem} (see, e.g., \cite[Theorem 1.25]{Shafarevich2013}) describes (or bounds) the dimension of algebraic varieties.
It arose in the context of the broader development of ``dimension theory'' in commutative algebra and algebraic geometry, which sought to formalize intuitive notions of dimension in the context of geometric objects.
The theorem lies at the basis of many dimension-counting arguments and is an everyday tool for algebraic geometers \cite[Section~12.5]{Clader2025}.

\begin{tool}[Fiber dimension; \Cref{thm:fiber_dimension}]
Let $\phi \colon X \to Y$ be a morphism of irreducible varieties.
\begin{enumerate}
\item For every $\vect{q} \in \phi(X)$, one has $\dim \phi^{-1}(\vect{q}) \geq \dim X - \dim \overline{\phi(X)}$.
\item There is a dense open subset $U \subseteq \overline{\phi(X)}$ such that $\phi$ has reduced fibers at $\vect{q} \in U$ and equality holds in the first statement.
\end{enumerate}
\end{tool}

The second tool, the \emph{Thom--Porteous formula} (see, e.g., \cite[Theorem 12.4]{EisenbudHarris2016}), also referred to as the Giambelli--Thom--Porteous formula, started with \cite{Giambelli1907}, and culminated in its modern form with \cite{Porteous1971}.
It gives an expression for the fundamental class of a degeneracy locus of a morphism of vector bundles in terms of Chern classes---we will clarify each part of this statement.

For our purposes, the vector bundle morphism can be viewed concretely as a parameter-dependent matrix (cf.~an MEP formulation): for each set of parameter values, the matrix defines a linear map between vector spaces of fixed dimension.
In this setting, the degeneracy locus of the morphism corresponds to the determinantal variety associated with this matrix, that is, the set of parameter values for which the matrix drops rank (cf.~the so-called eigentuples of the MEP). 
Such determinantal varieties are often singular, and their geometric structure can be difficult to study directly from their defining equations. 
Instead, the Thom--Porteous formula uses only structural information about the parameter-dependent matrix to compute basic invariants of the corresponding determinantal variety, providing a quantitative measure of \enquote{how large} and \enquote{how complicated} the variety is. 
In the special case where the determinantal variety is known to be zero-dimensional (i.e., it consists of finitely many points), the formula provides a way to count those points.

\begin{example}[Polynomial eigenvalue problems]\label{ex:PEP}
    For square one-parameter MEPs $\matx{M}(\lambda)\vect{v} = \vect{0}$, better known as polynomial eigenvalue problems, the determinantal variety is exactly the set of roots of the determinant $\det(\matx{M}(\lambda))$. 
    As the roots of a univariate polynomial, the eigenvalues are a finite set, whose count with multiplicity is determined by the degree of the determinant.
\end{example}

The Thom--Porteous formula gives a way to count solutions beyond polynomial eigenvalue problems.
Note that we only state the corank-one case; the general Thom--Porteous formula applies to arbitrary rank conditions, and for more general matrices (e.g., $e>f$ and/or mixed degrees, assuming their minors are still homogeneous).

\begin{tool}[Thom--Porteous formula (corank one); \Cref{thm:TP}]
Let $0 < e \leq f$ be integers, and let $\matx{M}(\vect{x}) \in \Cset[\vect{x}]^{f\times e}$ be a matrix of homogeneous polynomials in $\vect{x} \in \Pset^n$, whose columns have constant degree $d_1,\dots,d_e$. 
Assume that each component $Z \subseteq D_{e-1}(\matx{M})$ of the determinantal variety of rank $e-1$ has the expected dimension $\dim Z = n - (f-e+1) \geq 0$.
Then
\[
\deg D_{e-1}(\matx{M}) = T^{f-e+1}\text{-coefficient of }\frac{1}{(1-d_1T)\dotsm (1-d_eT)} \geq 0.
\]
\end{tool}

\subsection{Finiteness of the solution set}\label{sec:finiteness}

Analysis of the defining equations of $\mathcal{W}$ leads to the following structural description.

\begin{theorem}[Dimension count]\label{thm:V_sys_finite}
The Walsh variety $\mathcal{W}$ is an irreducible ideal-theoretic complete intersection in $\mathit{A^\circ\! BCDG}$ of dimension $2N+2$.
\end{theorem}

In simple terms, this means that the $N+n$ defining equations cut down the $(3N+n+2)$-dimensional ambient space $\mathit{A^\circ\! BCDG}$ to the expected dimension $2N+2$, and do so while intersecting transversally almost everywhere.
In order to avoid interrupting the main line of argument, we defer the proof to \Cref{sec:proof}.

In practice, we treat $(c,d)$ as the input data and $(a,b,g)$ as the unknowns.
For a specific problem instance $(c_0,d_0) \in CD$, the critical points are thus given by the intersection
\[
\mathcal{W}_{c_0,d_0} \coloneq \mathcal{W} \, \cap \, \{c = c_0, \,d = d_0\}.
\]  
Its dimension follows from the fiber dimension theorem. 

\begin{corollary}[Finiteness and reducedness of the critical points]\label{cor:finiteness}
Consider the projection 
\begin{equation}
\tau \, \colon \, \mathcal{W} \to \mathit{CD}, \, (a,b,c,d,g) \to (c,d). \label{eq:tau}
\end{equation}
For almost all $(c,d) \in CD$, the fiber $\tau^{-1}(c,d) = \mathcal{W}_{c,d}$ is empty or reduced of dimension one.
If it is non-empty, it consists of tuples $(a,b,c,d,g)$ satisfying
\begin{equation}
\deg a = n, \qquad a  \text{ has simple roots}, \qquad \gcd(a,b)=1. \label{eq:corFinit}
\end{equation}
In particular, for generic $H \in \stabModels{N}$, the set of complex solutions $(a,b,g)$ to the critical equations \eqref{eq:criticalEqsMatrix} is finite up to normalization of $a$ (e.g., fixing $a_n = 1$) and consists of simple solutions without multiplicities.
\end{corollary}

\begin{proof}
We distinguish two cases. 
When $\tau$ is not dominant, there is a dense open subset of points $(c,d) \in \mathit{CD}$ lying outside the image of $\tau$, so that for generic choice of $(c,d)$ we have $\tau^{-1}(c,d) = \emptyset$. 

Alternatively, if $\tau$ is dominant, i.e., $\overline{\tau(\mathcal{W})}= \mathit{CD}$, most fibers are not empty. 
In this case, the fiber dimension theorem (cf.~\Cref{thm:fiber_dimension}) implies that generic fibers are reduced of dimension
\[
\dim \tau^{-1}(c,d) = \dim \mathcal{W} - \dim \mathit{CD} = (2N+2)-(2N+1) = 1.
\]

The negations of the three additional properties claimed in the corollary are enforced by polynomial equations in the coefficients of $(a,b,g)$: equating the coefficient $a_n$, the discriminant $\Delta(a)$, or Sylvester's resultant $\operatorname{res}_s(a,b)$, respectively, to zero.
Since $\mathcal{W}$ is irreducible and none of these conditions is trivial\footnote{
One can easily construct a point $(a,b,c,d,g) \in \mathcal{W}$ that satisfies all of the conditions in \eqref{eq:corFinit}. For example, the point $a = (s+1)(s+2)\cdots(s+n)$ and $(b,c,d,g) = (1, 0, 0, 0)$.} on $\mathcal{W}$, they form a smaller-dimensional subset; of dimension at most $2N+1$. 
Let 
\[
\mathcal{V} \coloneq \mathcal{W} \cap \left( \{a_n = 0\} \cup \{\Delta(a) = 0\} \cup \{\operatorname{res}_s(a,b) = 0\} \right)
\]
and consider $\tau_{\mathcal{V}} \coloneq \tau|_{\mathcal{V}} \colon \mathcal{V} \to CD \colon (a,b,c,d,g) \mapsto (c,d)$. 
Since the solutions to the defining equations of $\mathcal{V}$ are invariant w.r.t.~the transformation $(a,b,c,d,g) \to (\lambda a,\lambda b, c, d, \lambda^{-1} g)$ for any $\lambda \in \Cset \setminus \{0\}$, all non-empty fibers $\tau_{\mathcal{V}}^{-1}(c,d)$ have dimension at least one.
Applying the fiber dimension theorem to every component of $\mathcal{V}$, one sees that $\dim \overline{\tau(\mathcal{V})} \le 2N < \dim \mathit{CD}$, hence $CD \setminus \tau(\mathcal{V})$ contains a dense open set.
We can conclude that for a generic choice of $(c,d)$, the fiber $\tau^{-1}(c,d)$ does not intersect $\mathcal{V}$ and thus $a,b,g$ satisfy the conditions in \eqref{eq:corFinit}.

Finally, since a generic fiber $\tau^{-1}(c,d)$ is one-dimensional and the condition $a_n \neq 0$ holds, we can normalize by setting $a_n=1$. 
This reduces the dimension of the generic fiber to zero, which completes the proof.
\end{proof}

\begin{remark}\label{rem:existence}
Note that this corollary does not claim the general fiber to be non-empty; though this will follow from \Cref{thm:numberSolutions} below.
In fact, it is known that there always exists a stable global minimizer of full McMillan degree exactly $n$ \cite[Theorem 1 \& Proposition 1]{Baratchart1986}. 
Baratchart works in the discrete-time setting, which is isometrically equivalent to the continuous-time setting (cf.~\Cref{sec:discretetime} for details).
\end{remark}

\subsection{Multiparameter eigenvalue problem}\label{sec:MEP}

The partially linear structure of the critical equations in \eqref{eq:criticalEqsMatrix} can be exploited to reformulate the critical equations as a multiparameter eigenvalue problem~\cite{Lagauw2023}. 
Let 
\[
\matx{M}(a,c,d) \coloneq \bmat{ \Toeplitz{a}{N}\tran \vect{d} & \Toeplitz{c}{n}\tran & \Toeplitz{\tilde{a}}{N}\tran\Toeplitz{\tilde{a}}{N-n}\tran} \in \Cset[\vect{a},\vect{c},\vect{d}]^{(N+n) \times (N+1)},
\]
so that the critical equations correspond to the following matrix-vector product:
\begin{equation}
   \Toeplitz{a}{N}\tran \vect{d} - \Toeplitz{c}{n}\tran \vect{b} - \Toeplitz{\tilde{a}}{N}\tran \Toeplitz{\tilde{a}}{N-n}\tran \vect{g} = \vect{0} \iff \matx{M}(a,c,d) \pmat{-1 \\ \vect{b} \\ \vect{g}} = \vect{0}.
    \label{eq:MEP}
\end{equation}
It follows that for each $(a,b,c,d,g) \in \mathcal{W}$, the vector $\Toeplitz{a}{N}\tran \vect{d}$ lies in the column space of  
\[
    \matx{B}(a,c) \coloneq \bmat{ \Toeplitz{c}{n}\tran & \Toeplitz{\tilde{a}}{N}\tran \Toeplitz{\tilde{a}}{N-n}\tran }, \label{eq:MORbmatrix}
\]
and hence $\matx{M}(a,c,d)$ fails to have full column rank at such points.

Fix input data $(c,d) \in CD$ and consider the determinantal variety determined by the maximal minors of $\matx{M}(a,c,d)$,
\begin{equation}
\mathcal{D} \coloneq D_N(\matx{M}(a,c,d)) = \Set{ [a] \mid \rank \matx{M}(a,c,d) \leq N } \subseteq \PP(A). \label{eq:detvar}
\end{equation}
Projectivizing is natural since all entries of $\matx{M}(a,c,d)$ are homogeneous in $a$ and this removes the scaling freedom in $a$. 
The points of $\mathcal{D}$ are the eigentuples of the MEP and can be computed using multiparameter eigenvalue solvers (e.g.,~\cite{Vermeersch2023,Hochstenbach2024,vermeersch2025a}).

It remains to show that the finiteness and reducedness of the generic fibers of $\mathcal{W}$ (\Cref{thm:V_sys_finite} and~\Cref{cor:finiteness}) carry over to the determinantal variety $\mathcal{D}$.
Define
\begin{equation}
    \phi \colon \mathit{A^\circ\!BG} \to \PP(A), \qquad  (a,b,g) \mapsto [a]. \label{eq:restriction}
\end{equation}
It is immediate from \eqref{eq:MEP} that $\phi(\mathcal{W}_{c,d}) \subseteq \mathcal{D}$.
However, $\mathcal{D}$ may contain additional components, corresponding to the closed subset
\[
\mathcal{B} \coloneq \Set{ [a] | \rank \matx{B}(a,c) < N } \subseteq \PP(A),
\]
which we refer to as the \emph{bad locus}.
This locus consists of points for which there exists a non-trivial vector in the null space of $\matx{M}(a,c,d)$ whose first entry vanishes. 
The corresponding solutions $(a,b,c,d,g)$ do \emph{not} satisfy the critical equations \eqref{eq:criticalEqsMatrix}, but instead solve the alternative system 
\[
\Toeplitz{c}{n}\tran \vect{b} + \Toeplitz{\tilde{a}}{N}\tran \Toeplitz{\tilde{a}}{N-n}\tran \vect{g} = \vect{0}.
\]
Consequently, we get 
\begin{equation}\label{eq:D_decomposition}
\mathcal{D} = \phi(\mathcal{W}_{c,d}) \cup \mathcal{B}.
\end{equation}
The following result shows that the bad locus $\mathcal{B}$ is in fact empty whenever $c$ has $N$ distinct roots. 
Under this assumption, $\mathcal{D}$ captures precisely the image of the critical points.


\begin{lemma}[Bad locus]\label{lem:MORbmatrix}
    Consider $(a,c) \in \mathit{A^\circ\!C^\circ}$ and let $\tilde a^\h$ and $c^\h$ denote the homogenization of the polynomials $\tilde a$ and $c$, respectively.
    Then,
    \begin{align}
        \nullity{\matx{B}(a, c)} &= \max \Big\{ 0, \; \deg \gcd(c, {\tilde a}^2) + \min \big\{N - \deg c - n, \; n - 2 \cdot \deg a \big\} \Big\} \nonumber \\
        &= \max \Big\{ 0, \; \deg \gcd(c^\h, (\tilde a^\h)^2) - n \Big\}. \label{eq:nullityB} 
    \end{align} 
    In particular, if $c \in C$ has $N$ distinct roots, then $\rank \matx{B}(a, c) = N$ for all $a \in \mathit{A^\circ}$.
    In that case, $\mathcal{B} = \emptyset$ and $\rank \matx{M}(a,c,d) \in \{N,N+1\}$.
\end{lemma}

The proof uses similar ideas to Bézout's lemma on polynomial greatest common divisors.

\begin{proof}
    Identifying $\vect{u} \in \Cset^n$ and $\vect{v} \in \Cset^{N-n}$ with $u^\h$ and $v^\h$, homogeneous polynomials of degree $n-1$ and $N-n-1$, respectively, we have
    \[
        \matx{B}(a,c) \pmat{\vect{u} \\ \vect{v}} = \vect{0} \iff c^\h u^\h + (\tilde a^\h )^2 v^\h = 0.
    \]
    Let $f \coloneq \gcd(c^\h, (\tilde a^\h)^2)$ and write $c^\h = f c'$ and $(\tilde a^\h)^2 = f a'$, so that $\gcd(c',a')=1$.
    Dividing by $f$ gives $u c' + v a' = 0$. 
    Since $a'$ divides $c' u$ and is coprime to $c'$, it divides $u$. 
    Hence $u = a' w$ and $v = -c' w$ for a $w \in \Cset[s,t]$ of degree $(n-1) - (2n - \deg f) = \deg f - n - 1$. 
    Conversely, every such $w$ yields a null vector of $\matx{B}(a,c)$ via the linear map $w \mapsto (a' w, -c' w)$.
    This map is hence an isomorphism from $\Cset[s,t]_{\deg f - n - 1}$ onto the null space of $\matx{B}(a,c)$. 
    Since $\Cset[s,t]_k$ has dimension $k+1$ for $k \geq 0$ and $\Cset[s,t]_k = \{0\}$ for $k < 0$, the nullity equals $\max\{0, \deg f - n\}$.

    If $c$ has $N$ distinct roots $p_1, \dots, p_N$, then $c^\h = c_N \prod_{i=1}^{N} (s - p_i t)$. 
    Since $f$ divides $c^\h$, it is a product of distinct factors $s - p_i t$, each of which divides $(\tilde a^\h)^2$. 
    Each factor being prime, it also divides $\tilde a^\h$,
    it follows that $f$ divides $\tilde a^\h$. In particular $\deg f \le n$, and the nullity in~\eqref{eq:nullityB} vanishes.    
\end{proof}



By~\Cref{lem:MORbmatrix}, the bad locus is empty whenever $c$ has $N$ distinct roots, so that every point of $\mathcal{D}$ corresponds to a critical point.\footnote{In model reduction practice, where $(c,d)$ correspond to the input data, the requirement $c_N \neq 0$ is trivially satisfied.
Furthermore, the higher-order model $H$ often originates from a system identification experiment based on noisy measurements, so that the estimated $c$ almost always has simple roots.}
The following corollary makes this precise: under the assumption, the critical equations and the determinantal variety describe the same solution set, preserving its reducedness, so that counting points of $\mathcal{D}$ is equivalent to counting the critical points.

\begin{corollary}[Equivalence of the MEP]\label{cor:equivalenceMEP}
    Consider input data $(c,d) \in CD$ with $c$ having $N$ distinct roots. 
    Then, the projection $\mathcal{W}_{c,d} \to \mathit{A^\circ}$ is an isomorphism onto the affine determinantal variety $\widehat{\mathcal{D}} = \set{a \in A^\circ | [a] \in \mathcal{D}}$.
    Consequently, if $(c,d)$ are generic, then $\mathcal{D}$ is a finite set of reduced points; that is, the maximal minors of $\matx{M}(a,c,d)$ define the elements in a reduced way.
\end{corollary}

\begin{proof}
    Fix $(c,d) \in CD$ with $c$ having $N$ distinct roots.
    For any $(a,b,g) \in \mathcal{W}_{c,d}$, the defining equations of $\mathcal{W}_{c,d}$ imply that $(-1,b,g)$ is a nonzero null vector of $\matx{M}(a,c,d)$, so $\matx{M}(a,c,d)$ is rank-deficient.
    Conversely, by the assumption on $c$ and \Cref{lem:MORbmatrix}, we have for any $a \in \widehat{\mathcal{D}}$ that $\rank \matx{B}(a,c) = N$, and hence the null space of $\matx{M}(a,c,d)$ is one-dimensional containing vectors of the form $(-\kappa, \kappa \vect{b}, \kappa \vect{g})\tran$, with $\kappa \in \Cset \setminus \{0\}$.  
    Upon normalization ($\kappa = 1$), this yields a unique pair $(b,g)$ for which $(a,b,g) \in \mathcal{W}_{c,d}$. 
    Hence the inverse map $a \mapsto (a,b,g)$ is a morphism, given explicitly by Cramer's rule.
    This establishes the isomorphism.

    \Cref{cor:finiteness} further shows that, for generic $(c,d)$, the fiber $\mathcal{W}_{c,d}$ is finite up to scaling of $a$.
    Moreover, the same corollary asserts that the critical equations define the generic fibers in a reduced way.
    By \Cref{lem:reducedMinors}, this reducedness carries over to the determinantal description.
\end{proof}

In fact, this argument yields a proof of finiteness of the critical locus for \emph{arbitrary} higher-order models of McMillan degree exactly $N$ and distinct poles.
The following result sharpens \Cref{cor:finiteness}, and yields a proof independent of \cite[Theorem~3.2]{Hanzon2012}; note that it does not rely on \Cref{thm:V_sys_finite} and hence is rather direct.

\begin{corollary}\label{cor:finiteness_by_simple_roots}
If $c$ has distinct roots and $\gcd(c,d) = 1$, then $\mathcal{D}$ is a finite set.
\end{corollary}
\begin{proof}
Since $c$ has $N$ distinct roots, \Cref{cor:equivalenceMEP} provides the inverse morphism $\widehat{\mathcal D} \to \mathcal W_{c,d}$, $a \mapsto (a,b,g)$, defined on all of $\widehat{\mathcal D}$. For $0 \le i \le n$ and $0 \le j \le N-n-1$, consider the composite map
\[
  \varphi_{ij}\colon \widehat{\mathcal D} \to \mathcal W_{c,d} \to \Cset, \qquad a \mapsto (a,b,g) \mapsto a_i \cdot g_j .
\]
Scaling $a$ by $\lambda\in \Cset \setminus \{0\}$ scales $g$ by $\lambda^{-1}$, that is, $\varphi_{ij}(\lambda a) = \varphi_{ij}(a)$, so $\varphi_{ij}$ is a well-defined regular function on $\mathcal{D}$. 
Since $\mathcal{D} \subseteq \PP(A)$ is a projective variety, \Cref{thm:constant} indicates that $\varphi_{ij}$ must be \emph{constant} on each irreducible component of $\mathcal{D}$. 

If $[a] \in \mathcal{D}$, then $g$ is not identically zero, since by \eqref{eq:lofs} we would otherwise have $ad=bc$, so $c \mid a$ by $\gcd(c,d)=1$, contradicting $\deg c = N > n \geq \deg a$ (and $a\neq 0$). 
In particular, for any $[a] \in \mathcal D$ there exist $i,j$ such that $\varphi_{i,j}(a) \neq 0$. 
In a neighborhood of $[a]$ in $\mathcal{D}$, the coordinate ratio
\[
[a] \mapsto \frac{a_k}{a_i} =  \frac{a_k\cdot g_j}{a_i \cdot g_j} = \frac{\varphi_{kj}(a)}{\varphi_{ij}(a)}, \qquad k=0,\dots,n
\]
is thus constant, so $[a]$ is isolated in $\mathcal{D}$. This shows that $\mathcal{D} \subseteq \PP(A)$ is discrete, hence finite.
\end{proof}

\subsection{Degree of the solution set}\label{sec:count}

To obtain an explicit count of the critical points, it remains to apply the Thom--Porteous formula (\Cref{thm:TP}) to the determinantal variety $\mathcal{D}$.

\begin{theorem}[Solution count]\label{thm:numberSolutions} 
Consider input data $(c,d) \in CD$ such that $c$ has $N$ distinct roots and $\gcd(c,d)=1$, and let $m$ denote the number of critical points, that is, the number of complex solutions to the critical equations \eqref{eq:criticalEqsMatrix} up to scaling.
We have
\begin{equation}
    m = \# \mathcal{D} \leq \deg \mathcal{D} = \sum_{j=0}^n \binom{N-n-1+j}{j} 2^j \label{eq:numSolutions}
\end{equation} 
and equality holds for generic $(c,d)$.
\end{theorem}
\begin{proof}
The expected dimension of $\mathcal{D}$ is $n - (N+1-N)(N+n-N) = 0$.
By \Cref{cor:finiteness_by_simple_roots}, $\mathcal{D}$ is indeed finite, so that each component has the expected dimension zero, and we can apply the Thom--Porteous formula for the corank one case (cf.~\Cref{thm:TP}). 
The coefficient of $T^n$ in the series expansion of $(1 - T)^{-1}(1 - 2T)^{-(N - n)}$ can be determined by expressing each factor as a power series. Specifically,
\[
    \frac{1}{1 - T} = \sum_{j = 0}^{\infty} T^j,
    \qquad
    \frac{1}{(1 - 2T)^{N - n}} = \sum_{j = 0}^{\infty} \binom{N - n + j - 1}{j} (2T)^j,
\]
such that by the Cauchy product
\[
    \frac{1}{(1 - T)(1 - 2T)^{N - n}} = \sum_{k = 0}^{\infty} \left( \sum_{j = 0}^{k} \binom{N - n + j - 1}{j} 2^j \right) T^k. 
\] 
In this series, the coefficient of $T^n$ corresponds to
\[
    \deg \mathcal{D} = \sum_{j=0}^n \binom{N-n-1+j}{j} 2^j.
\]

By~\Cref{cor:equivalenceMEP}, $\mathcal{D}$ is in bijection with the critical points, so that $m = \#\mathcal{D} \leq \deg \mathcal{D}$.
The same corollary shows that for generic $(c,d)$, the maximal minors define $\mathcal{D}$ in a reduced way, so that $\#\mathcal{D} = \deg \mathcal{D}$.
\end{proof}

Several values of the formula from \Cref{thm:numberSolutions} are depicted in \Cref{tab:numSols}. 
If we consider the co-order one case $n = N-1$, then the formula \eqref{eq:numSolutions} reduces to $\sum_{j=0}^n 2^j = 2^{n+1}-1$, which corresponds exactly to the upper bound obtained in \cite{Hanzon2007}.

\begin{table}[h]
    \centering
    \caption{Evaluation of the formula from \Cref{thm:numberSolutions} for various problem sizes $(N,n)$.
    Under the assumption that $c,d$ are generic, each entry gives the number of complex solutions to the critical equations of the model reduction problem.
    }
    \vskip1em
    \begin{tabular}{c|rrrrrrrrrrr}
    \diagbox[height=1.4\line]{$N$}{$n$} &    1  &     2  &      3  &      4  &       5  &       6  &       7  &       8  &       9  &       10  \\
    \hline
    2  &   3  &        &         &         &          &          &          &          &          &                   \\
    3  &   5  &    7   &         &         &          &          &          &          &          &                   \\
    4  &   7  &   17   &    15   &         &          &          &          &          &          &                   \\
    5  &   9  &   31   &    49   &    31   &          &          &          &          &          &                   \\
    6  &  11  &   49   &   111   &   129   &     63   &          &          &          &          &                   \\
    7  &  13  &   71   &   209   &   351   &    321   &    127   &          &          &          &                   \\
    8  &  15  &   97   &   351   &   769   &   1023   &    769   &    255   &          &          &                   \\
    9  &  17  &  127   &   545   &  1471   &   2561   &   2815   &   1793   &    511   &          &                   \\
    10 &  19  &  161   &   799   &  2561   &   5503   &   7937   &   7423   &   4097   &   1023   &                   \\
    11 &  21  &  199   &  1121   &  4159   &  10625   &  18943   &  23297   &  18943   &   9217   &   2047            
\end{tabular}
 \label{tab:numSols}
\end{table}

To put these numbers in perspective, we compare~\eqref{eq:numSolutions} with the alternative bound provided by B\'ezout's theorem. 
Applying that result to the critical equations~\eqref{eq:criticalEqsMatrix}, in the chart $a_n = 1$, gives an upper bound of $2 \cdot 3^{N+n-2}$ complex solutions, since one equation is linear, one is quadratic, and the remaining $N+n-2$ are cubic. 
For $(N,n) = (6,2)$, the Bézout and Thom--Porteous numbers are $1458$ and $49$, respectively; for $(N,n) = (10,5)$, they are $3\,188\,646$ and $5503$.

\begin{remark}
In contrast to the number of complex critical points, the number of real-valued and Hurwitz stable solutions---corresponding to the genuine stationary points of the optimization problem, e.g., local minima, saddle points---depends on the specific input data $c,d$ and cannot, in general, be determined solely from the problem size $(N,n)$.
\end{remark}

\subsection{Irreducibility of the Walsh variety}\label{sec:proof}

We complete this section by proving \Cref{thm:V_sys_finite}. 
We first point out a useful structure in the defining equations of the Walsh variety $\mathcal{W}$.
Consider the projection 
\begin{equation}
\pi \colon \mathcal{W} \to \mathit{A^{\circ}\!BG},\; (a,b,c,d,g) \to (a,b,g).\label{eq:pi}
\end{equation}
For fixed $(a,b,g)$, the right-hand side of the defining equation of $\mathcal{W}$ is constant, so the unknowns $(c,d)$ are constrained by an affine system of linear equations:
\begin{equation}
    \bmat{\Toeplitz{a}{N}\tran & -\Toeplitz{b}{N+1}\tran}\pmat{\vect{d} \\ \vect{c}} = \Toeplitz{\tilde{a}}{N}\tran \Toeplitz{\tilde{a}}{N-n}\tran \vect{g}. \label{eq:lofsMatrix}
\end{equation} 
Consequently, each fiber of $\pi$ is an affine linear subspace:
\begin{equation} \label{eq:affinebundle}
\pi^{-1}(a,b,g) = \{(a,b,g)\} \times \Set{(c,d) | \eqref{eq:lofsMatrix} \text{ is satisfied}}.
\end{equation}

\begin{lemma}\label{lem:divisibilityConstraint}
For given $(a,b,g) \in \mathit{A^{\circ}\!BG}$, the fiber $\pi^{-1}(a,b,g)$ is non-empty if and only if $\gcd(a,b)$ divides $\tilde{a}^2 g$.
In that case, we have
\begin{align*}
\dim \pi^{-1}(a,b,g) &=  N + \deg \gcd(a,b) - \max \{ \deg b, \deg a  - 1 \} \\
 &= N + \deg \gcd(a^{\mathrm h},b^{\mathrm h}) - (n-1),
\end{align*}
where we use the convention $\deg 0 = -\infty$ and, as usual, $\gcd(a,0) = a$.
\end{lemma}

The notation and construction in the proof is similar to \Cref{lem:MORbmatrix}.

\begin{proof}
Let $f \coloneq \gcd(a,b)$ and $k \coloneq \deg f$, and write $a' \coloneq  a/f$ and $b' \coloneq b/f$ so that $\gcd(a',b') = 1$. 
It is immediate from
\[
 f (a' d - b' c) = \tilde{a}^2 g 
\]
that $f \mid \tilde{a}^2 g $ is necessary for the fiber \eqref{eq:affinebundle} to be non-empty. 

We show that it is also sufficient. 
Let $\matx{S}(a,b) \coloneq \bmat{\Toeplitz{a}{N}\tran & -\Toeplitz{b}{N+1}\tran} \in \Cset^{(N+n) \times (2N+1)}$. 
An element $(\vect{u}\tran, \vect{v}\tran)\tran \in \Cset^{2N+1}$ in the null space of $\matx{S}(a,b)$ corresponds to polynomials $u \in \CC[s]_{\leq N-1}$ and $v \in \CC[s]_{\leq N}$ satisfying $a u - b v = 0$.
Removal of the common factor gives $a' u = b' v $, and since $\gcd(a',b') = 1$, this holds if and only if there exists a polynomial $r$ such that $u = b' r$ and $v = a' r$.
Taking the degrees of $u$ and $v$ into account gives
\[
\deg r \leq \rho \coloneq \min \big\{\underbrace{(N - 1)}_{\deg u } - \underbrace{(\deg b  - k)}_{\deg b' }, \; \underbrace{(N)}_{\deg v } - \underbrace{(\deg a  - k)}_{\deg a' }\big\},
\]
so that $r \mapsto (b'r, a'r)$ is an isomorphism from $\Cset[s]_{\leq\rho}$ onto the null space of $\matx S(a,b)$.
Hence we obtain the first expression in the lemma:
\begin{equation*}
\nullity{\matx{S}(a,b)} = \rho + 1 = N + k - \max \{ \deg b, \deg a  - 1 \}.
\end{equation*}
Since $\gcd(a^\h, b^\h) = t^{\min\{n-\deg a, n-1-\deg b\}} f^\h$ (where $f^\h$ is homogenized to degree $\deg f$), we have 
\[
\deg\gcd(a^\h, b^\h) = k + n - 1 - \max\{\deg b, \deg a - 1\},
\] which gives the second expression for the dimension in the lemma. 
Taking into account the size of $\matx S(a,b)$, we obtain\footnote{Note that \eqref{eq:rankS} generalizes the basic fact that the Sylvester matrix of $a$ and $b$ has nullity equal to $\deg\gcd(a^\h, b^\h)$.}
\begin{equation}
  \rank \matx S(a,b) = N + \max\{\deg b, \deg a - 1\} - k + 1 = N + n - \deg\gcd(a^\h, b^\h). \label{eq:rankS}
\end{equation}

The image of $\matx{S}(a,b)$ is the subspace of coefficient vectors corresponding to polynomials $p a - q b$ for suitable $p \in \CC[s]_{\leq N-1}$ and $q \in\CC[s]_{\leq N}$.
Each such polynomial is divisible by $f$ and has degree at most $N + \max\{\deg b, \deg a - 1\}$, so that
\[
\image \matx{S}(a,b) \subseteq \Set{ f \cdot h | h\in  \CC[s]_{\leq N+\max \{ \deg b, \; \deg a  - 1 \} - k} }.
\] 
By \eqref{eq:rankS}, both sides have dimension $N+ \max \{ \deg b, \; \deg a  - 1 \} - k + 1$, so that this inclusion is in fact an equality.
Moreover, since $\deg a \le n$, we have $\deg(\tilde a^2 g) \le 2\deg a + N - n - 1 \le N + \deg a - 1$. 
Hence, if $f$ divides $\tilde a^2 g$, then $\tilde a^2 g$ lies in the image of $\matx S(a,b)$, which shows that the divisibility constraint is also sufficient. 
If it is satisfied, the affine space of solutions $\pi^{-1}(a,b,g)$ has dimension $\nullity{\matx S(a,b)}$.
\end{proof}

Now that we have established \Cref{lem:divisibilityConstraint}, we can present the proof of \Cref{thm:V_sys_finite}.
It proceeds in the following steps.
First, we stratify $\mathcal{W}$ by the degree of the greatest common divisor of $a$ and $b$, this allows us to determine the dimension of $\mathcal{W}$. 
Then, we show that $\mathcal{W}$ is a set-theoretic complete intersection and irreducible. Lastly, reducedness follows from a bundle argument, so $\mathcal{W}$ is an ideal-theoretic complete intersection.

\begin{proof}[Proof of \Cref{thm:V_sys_finite}]
For $0\leq k \leq n$, consider the following subsets of $\mathit{A^\circ\!BG}$:
\[
\mathcal{S}_{\geq k} \coloneq \Set{(a,b,g) | \deg \gcd(a^\h,b^\h) \geq k }, \qquad \mathcal{S}_k \coloneq \mathcal{S}_{\geq k} \setminus \mathcal{S}_{\geq k+1},
\]
\[
\mathcal{T}_k \coloneq \Set{(a,b,g) \in \mathcal{S}_k | \gcd(a,b) \text{ divides } \tilde{a}^2 g}.
\]
The set $\mathcal{S}_{\geq k}$ is defined by the vanishing of the $(2n-k)$-minors of the Sylvester matrix constructed from the polynomials $a^\h$ and $b^\h$.
On $\mathcal{S}_k$, we have $\rank \matx S(a,b) = N+n-k$ by~\eqref{eq:rankS}, and, by the proof of~\Cref{lem:divisibilityConstraint}, the divisibility constraint holds if and only if the coefficient vector of $\tilde a^2 g$ lies in the column space of $\matx S(a,b)$.
Within $\mathcal{S}_k$, the set $\mathcal{T}_k$ is thus defined by the vanishing of the $(N+n-k+1)$-minors of the matrix $\matx{S}(a,b)$ concatenated with the coefficient vector of $\tilde{a}^2g$.

Using the convention $\Cset[s]_{\le -1} \coloneq \{0\}$, the variety $\mathcal{S}_{\geq k}$ is the image of the polynomial map
\[
(\CC[s]_{\leq k} \setminus 0) \times (\CC[s]_{\leq n-k} \setminus 0) \times \CC[s]_{\leq n-k-1} \times G \to \mathit{A^\circ\!BG}, \qquad (f,u,v,g) \mapsto (fu,fv,g).
\]
The image of this map is exactly $S_{\ge k}$, which is therefore irreducible.

Given $(a,b,g) \in \mathcal{S}_{\geq k}$, there are only finitely many factorizations $a=fu$, $b=fv$ up to rescaling $f$ by $\lambda \in \Cset \setminus \{0\}$. 
Hence, the fibers of this map are $1$-dimensional, so that $\mathcal{S}_{\geq k}$ is an irreducible variety of dimension $N+n+1-k$ (\Cref{thm:fiber_dimension}).
Since $S_{\geq k+1}$ has smaller dimension, $S_k$ is a non-empty open subset of $S_{\geq k}$, and thus irreducible of the same dimension.

If $k=0$, then the divisibility constraint is trivial and $\mathcal{T}_0 = \mathcal{S}_0$, a non-empty open subset of $\mathit{A^{\circ}\!BG}$.
For $k\geq 1$, $\mathcal{T}_k$ is a proper closed subset of $\mathcal{S}_k$, since there are tuples $(a,b,g) \in \mathcal{S}_k$ that do \emph{not} satisfy the divisibility constraint, for example,
\[
a=(s-1)^n, \quad b=(s-1)^{k}(s-2)^{n-1-k}, \quad g=(s-3)^{N-n-1},
\]
for $1 \le k \le n-1$, and similarly for $k = n$.
Therefore $\dim \mathcal{T}_0  = N+n+1$ and $\dim \mathcal{T}_k \leq N+n-k$ for $k\geq 1$.
Consider again the map $\pi$ from \eqref{eq:pi}; we stratify $\mathcal{W}$ according to the degree of $\gcd(a^\h,b^\h)$:
\[
\mathcal{W} = \bigcup_{k=0}^{n} \mathcal{W}_k, \qquad \mathcal{W}_k \coloneq \Set{ (a,b,c,d,g)\in\mathcal{W} | \deg \gcd(a^\h,b^\h)=k} = \pi^{-1}(\mathcal{T}_k)
\]
where the last equality holds by~\Cref{lem:divisibilityConstraint}.
By said lemma, the fiber of $\pi$ over $\mathcal{T}_k$ are the solution sets of the linear system~\eqref{eq:lofsMatrix} of constant rank $N+n-k$.
Hence, $\pi \colon \mathcal{W}_k \to \mathcal{T}_k$ is an affine bundle of rank $N+k-(n-1)$, and
\[
\dim \mathcal{W}_k = \dim \mathcal{T}_k + \rank_{\mathcal{T}_k} \mathcal{W}_k
\leq (N+n-k+1) + (N+k-(n-1)) = 2N+2, 
\] 
where equality holds exactly for $k=0$. This establishes $\dim \mathcal{W} = \max_{k} \dim \mathcal{W}_k = 2N+2$.
Moreover, $\mathcal W_0$, as an affine bundle over the irreducible base $T_0$, is irreducible.


From this we can conclude that $\mathcal{W}$ is a set-theoretic complete intersection:
\[
\codim_{\mathit{A^\circ\! BCDG}} \mathcal{W} = \underbrace{(n+1)+n+(N+1)+N+(N-n)}_{\text{ambient dimension}} - \underbrace{2N+2}_{\dim \mathcal{W}} = \underbrace{N+n}_{\!\!\#\text{equations}\!\!}\!\!.
\]
Hence, by \Cref{thm:intersection_dim}, every irreducible component of $\mathcal{W}$ has dimension exactly $2N+2$. 

On the other hand, $\mathcal{W} = \bigcup_{k=0}^{n} \overline{\mathcal{W}_k}$ is a union of subvarieties, so none of the $\mathcal{W}_k$, $k\geq 1$, are (unions of) components of $\mathcal{W}$, as they have strictly smaller dimension. 
In this way, we conclude that $\mathcal{W} = \overline{\mathcal{W}_0}$ is irreducible.
Finally, we have established that $\mathcal{W}_0$ is a dense open subset of $\mathcal{W}$, and it is a bundle over the smooth variety $\mathcal{T}_0$, an open subset of $\Cset^{N+n+1}$.
As such, $\mathcal{W}$ is reduced (even smooth) along $\mathcal{W}_0$ and by the Unmixedness \Cref{thm:unmixedness}, the variety $\mathcal{W}$ is reduced everywhere; meaning that the critical equations define $\mathcal{W}$ in a reduced way.
\end{proof}

\section{Further considerations}\label{sec:otherConsid}
We present extensions and refinements of the results from \Cref{sec:finitAndDeg}.

\subsection{Full-order models with multiple poles}\label{sec:multiplepoles}
If $\deg c = N$ but $c$ has multiple roots, the bad locus $\mathcal{B}$ may be finite or even positive-dimensional. 
This is illustrated in the following example.

\begin{example}[Multiple poles]\label{H2MOR_TP:multiplePoles}
    Apply the affine chart $a_1 = 1$ and consider $(N,n) = (2,1)$ with $c = (s-p)^2$ for some $p \in \Cset$.
    In this case, $\matx{B}(a, c)$ drops rank if and only if $a = (s+p)$.
    Alternatively, let $(N,n) = (4,3)$ and $c = (s-p_1)^2 (s-p_2)^2$ for some $p_1,p_2 \in \Cset$. 
    Here, $\mathcal{B}$ consists of polynomials $a = (s+p_1)(s+p_2)(s-p_3)$, for arbitrary $p_3 \in \CC \cup \{\infty\}$. \qedex
\end{example}

When the bad locus is non-empty, it is natural to ask whether the images of the critical points can still be distinguished from it.
The following corollary shows that for reasonable input data---
higher-order models of McMillan degree exactly $N$---the bad locus never overlaps with the images of the critical points.

\begin{corollary}[Critical points avoid the bad locus]\label{cor:disjoint}
    Let $\phi$ be the projection from \eqref{eq:restriction} and consider input data $(c,d) \in CD$ with $\deg c = N$ and $\gcd(c,d) = 1$. 
    Then $\phi(\mathcal{W}_{c,d}) \cap \mathcal{B} = \emptyset$, so that so that every point of $\mathcal{D}$ is either the image of a critical point or an element of the bad locus, but not both.
    Moreover, for every $[a] \in \phi(\mathcal{W}_{c,d})$, the null space of $\matx{M}(a,c,d)$ is spanned by a single vector $(-1, b, g)$, so that the corresponding critical point $(a,b,g)$ is unique up to scaling.
    \end{corollary}
    
    \begin{proof}
    Let $(a,b,g) \in \mathcal{W}_{c,d}$ and $f \coloneq \gcd(c, \tilde a^2)$. 
    Since $ad = bc + \tilde a^2 g$, the polynomial $f$ divides $ad$. 
    By the assumption, $\gcd(c,d) = 1$, so that $f$ is coprime to $d$ and hence $f$ divides $a$. 
    Therefore $\deg f \le \deg a \le n$. 
    Since $\deg c = N$, $t$ does not divide $c^\h$, so that $\gcd(c^\h, (\tilde a^\h)^2)$ also has degree $\deg f \le n$. 
    This implies, by~\Cref{lem:MORbmatrix}, that $\matx{B}(a,c)$ has full column rank $N$, i.e., $[a] \notin \mathcal{B}$. Consequently, the null space of $\matx{M}(a,c,d) = \bmat{ \Toeplitz{a}{N}\tran \vect{d} & \matx{B}(a,c) }$ is at most
    one-dimensional, and any nonzero null vector has a nonzero first entry. 
    It therefore contains $(-1, b, g)$ and is spanned by it.
\end{proof}

    If $\mathcal{B}$ is non-empty but finite, part of the generic count for $\deg \mathcal{D}$ in~\Cref{thm:numberSolutions} is accounted for by points in the bad locus; subtracting their contribution decreases the count for the critical points accordingly. 
    In this case, the combination of~\Cref{cor:disjoint} and~\Cref{thm:numberSolutions} provides an upper bound on the number of critical points.

    For $[a] \in \mathcal{D}$, let $\operatorname{mult}_a I_{N+1}(\matx{M}(a,c,d))$ denote its multiplicity with respect to the ideal generated by the maximal minors of $\matx{M}(a,c,d)$, and let $\deg \mathcal{B}$ denote the degree of the ideal $I_N(\matx{B})$ generated by the maximal minors of $\matx{B}(a,c)$.

    \begin{corollary}[Adjusted count for a non-empty bad locus]\label{cor:multbound}
        Consider input data $(c,d) \in CD$ with $\deg c = N$ and $\gcd(c,d) = 1$, and assume that $\mathcal{D}$ is finite. 
        Then the number $m$ of critical points satisfies
        \begin{equation}
            m = \# \mathcal{D} - \# \mathcal{B} \leq \deg \mathcal{D} - \sum_{[a] \in \mathcal{B}} \operatorname{mult}_a I_{N+1}(\matx{M}(a,c,d)) \leq \deg \mathcal{D} - \deg \mathcal{B}, \label{eq:multbound}
        \end{equation}
        where $\deg \mathcal{D}$ is still given by \eqref{eq:numSolutions}. 
        Equality holds in the first inequality if and only if every point of $\mathcal{D} \setminus \mathcal{B}$ has multiplicity one.
    \end{corollary}
    \begin{proof}
        By the decomposition $\mathcal{D} = \phi(\mathcal{W}_{c,d}) \cup \mathcal{B}$ and~\Cref{cor:disjoint}, the critical points are in bijection with the points of $\mathcal{D} \setminus \mathcal{B}$, so that $m = \#\mathcal{D} - \#\mathcal{B}$. 
        Since $\mathcal{D}$ is finite, the Thom--Porteous formula applies and $\deg \mathcal{D}$ equals~\eqref{eq:numSolutions}. 
        Each point of $\mathcal{D}$ contributes its multiplicity, which is at least one, to $\deg \mathcal{D}$. 
        This gives the first inequality, with equality if and only if every point of $\mathcal{D} \setminus \mathcal{B}$ has multiplicity one. 
        The second inequality holds because $I_{N+1}(\matx{M}(a,c,d)) \subseteq I_N(\matx{B}(a,c))$, so that the multiplicity of every $[a] \in \mathcal{B}$ with respect to the maximal minors of $\matx{M}(a,c,d)$ is bounded below by the multiplicity of the point with respect to the maximal minors of $\matx{B}(a,c)$.
    \end{proof}

In particular, for every given model $H \in \mathcal{M}_N$ for which $\mathcal{D}$ is finite, the count~\eqref{eq:numSolutions} remains an upper bound on the number of critical points, also when $H$ has multiple poles. One can deduce from \Cref{lem:MORbmatrix} that $\mathcal{B}$ is finite if and only if $c$ has at most one multiple root, which may be of multiplicity at most $3$. We analyze the example of one double root below.

\begin{proposition}\label{H2MOR_TP:multiplerootscount}
Let $(c,d) \in CD$, where $c$ has exactly one double root and $N-2$ simple roots. 
Then $\# \mathcal{B} = \binom{N-2}{n-1}$.
If, moreover, $\gcd(c,d)=1$ and $\mathcal{D}$ is finite, then the number $m$ of critical points satisfies
\[
m   = \# \mathcal{D} - \binom{N-2}{n-1} \leq \sum_{j=0}^n \binom{N-n-1+j}{j} 2^j - \binom{N-2}{n-1}.
\]
\end{proposition}

\begin{proof}
Write $c = (s-p_1)^2 \prod_{i=2}^{N-1}(s-p_i)$ with $p_1,\dots,p_{N-1}$ distinct. 
The set $\mathcal{D}$ can be split into critical points and bad locus points. 
From \Cref{lem:MORbmatrix} we deduce that $\rank \matx{B}(a,c) < N$ happens if and only if $\tilde{a} = (s-p_1) \prod_{i \in I} (s-p_i)$, where $I \subseteq \{2,\dots,N-1\}$ is a selection of $n-1$ roots of $c$. 
Thus $\mathcal{B}$ consists of exactly $\binom{N-2}{n-1}$ points, disregarding potential multiplicity. 
The result is then an application of \Cref{cor:multbound}.
\end{proof}

\subsection{Discrete-time models}\label{sec:discretetime}

All results presented in this paper have a direct counterpart in the discrete-time setting. 
They are obtained via identical derivations but adapted notation, or by direct transformation of our results via the isomorphism on spaces $\mathcal{H}_2(\mathbb{C}^{+})$ and $\mathcal{H}_2(\mathbb{D})$, the Hardy spaces on the half plane and the unit disc, respectively~\cite[Theorem 1.2.5]{Partington2004}.
A treatment of this discrete-time case can be found in \cite[Chapter~5]{Leonie2025}.

In the discrete-time case, the transfer functions under consideration are strictly proper rational functions of the complex variable $z$, analytic outside the closed unit disk. 
Given such an $N$th order discrete-time model $H \in \stabDiscModels{N}$ and a desired reduced order $1 \leq n < N$, the $\Hardy$-norm model reduction problem searches for a solution
\[
\widehat{H}(z) \in \operatorname*{argmin}_{\widehat{H} \in \stabDiscModels{n}} \; J_{\mathrm{d}}^2,
\]
where the subscript $\mathrm{d}$ indicates the discrete-time domain. 
For example,
\[
J_{\mathrm{d}}^2 = \norm{ H(z) - \widehat{H}(z)}^2_{\mathcal{H}_2} = \frac{1}{2\pi} \int_{0}^{2\pi}
\left| H (e^{\mathrm{j}\omega}) - \widehat{H}(e^{\mathrm{j}\omega}) \right|^2 \, \diff \omega = \sum_{k = 1}^\infty \left( h_k - \widehat{h}_k \right)^2.
\]
The samples $\{h_k\}_{k=1, \dots,\infty}$ and $\{\widehat{h}_k\}_{k=1, \dots,\infty}$ are the impulse responses of $H$ and $\widehat{H}$, respectively.

In this setting, the interpolatory conditions in \Cref{thm:interpolation,H2MOR_thm:walsh} change: instead of reflecting the poles across the imaginary axis, one now inverts them with respect to the unit circle, so that the interpolation points $- \hat p_i$ are replaced by $1/ \hat p_i$.
Accordingly, our notation $\tilde{a} = \sum_{i=0}^n a_i(-s)^i$ used in \eqref{eq:lofs} must be adapted to the discrete-time setting and becomes $\tilde{a}(s) = z^{-n}a(1/z) = \sum_{i=0}^n a_i z^{n-i}$.

\begin{theorem}[Discrete-time count]
    With adapted notation, all results in \Cref{sec:finitAndDeg,sec:multiplepoles}, including \Cref{cor:finiteness,lem:MORbmatrix,cor:finiteness_by_simple_roots,thm:numberSolutions,cor:disjoint,cor:multbound}, hold verbatim for the discrete-time model order reduction problem.
\end{theorem}

\section{Examples: the good, the bad and the degenerate}\label{sec:Numerical}

We present some (numerical) examples to illustrate our results. 
The code is available in a public repository \cite{h2morcode}.
Throughout, numerical values are rounded for displaying purposes.

We start with \Cref{ex:easy,H2MOR_ex:didactic}, illustrating the behavior for generic input data; the remaining examples consider special behavior.
In \Cref{ex:mulitpole_higher_mult,ex:badlocusmult} the input data is such that the bad locus is non-empty, whose presence turns the exact count into an upper bound (\Cref{sec:multiplepoles}). 
\Cref{ex:critMult} considers special input data for which a critical point appears with multiplicity, contrary to the generic case characterized in \Cref{cor:finiteness}. 
\Cref{exmp:vandooren,exmp:enforcedMult} investigate the subtle differences between the different optimality conditions from \Cref{thm:interpolation,H2MOR_thm:walsh} for approximants with multiple poles.
Finally, \Cref{exmp:enforcedMult} and \Cref{rem:degen} briefly touch upon the phenomenon of \emph{degenerate} stationary points, which are reduced models at which the gradient of $J^2$ with respect to the coefficients of $(a,b)$ vanishes, but $\gcd(a,d) \neq 1$ so that $\widehat{H} \notin \mathcal{M}_n$.

\begin{example}[Continuous time]\label{ex:easy}
    We revisit the example of \cite[Section 5]{Agudelo2021} and \cite[Example 1]{Lagauw2023}, which seeks the optimal first-order approximation ($n=1$) of the third-order model ($N=3$)
    \[
    H(s) = \frac{d}{c} = \frac{s^2 + 9s-10}{s^3+12s^2+49s+78} = \frac{(s+10)(s-1)}{(s+6)(s^2+6s+13)} \in \mathcal{M}_3,
    \]
    whose poles are $s \in \{-6,-3\pm 2\mathrm{j}\}$.
    For this setup, the matrix $\matx{M}(a,c,d)$ corresponds to
    \[
    \matx{M}(a,c,d) = \bmat{
        -10a_0 & 78 & a_0^2 & 0 \\
        9a_0 -10a_1  & 49 & -2a_0a_1 & a_0^2 \\ 
        a_0 + 9a_1 & 12 & a_1^2 & -2a_0a_1 \\ 
        a_1 & 1 & 0 & a_1^2
    }.
    \]
    The determinantal variety $\mathcal{D}$, defined by the determinant of $\matx{M}(a,c,d)$, consists of five reduced points $[a] = (a_0 : a_1) \in \Pset^1$:
    \[
    \mathcal{D} = \Set{(9.6796 : 1),(0.2671 : 1),(-16.6189 : 1),(-4.1639 + 0.9027 \mathrm{j} : 1),(-4.1639 - 0.9027\mathrm{j} : 1)}.
    \]
    Since $c$ has $N$ distinct roots, the bad locus $\mathcal{B}$ is empty (cf.~\Cref{lem:MORbmatrix}).
    The number of obtained critical points is in agreement with the count from \Cref{thm:numberSolutions}.
    The two real points in $\mathcal{D}$ with $a_0/a_1 > 0$ correspond to models 
    \[
    \widehat{H}_1 = \frac{1.279}{s+9.679}, \qquad \widehat{H}_2 = \frac{-0.043}{s+0.267}
    \]
    in $\mathcal{M}_1$, both of which are local minimizers of the $\mathcal{H}_2$-error; $\widehat{H}_1$ being the global minimizer.\qedex
\end{example}

\begin{figure}
    \centering
    \includegraphics[width=\textwidth]{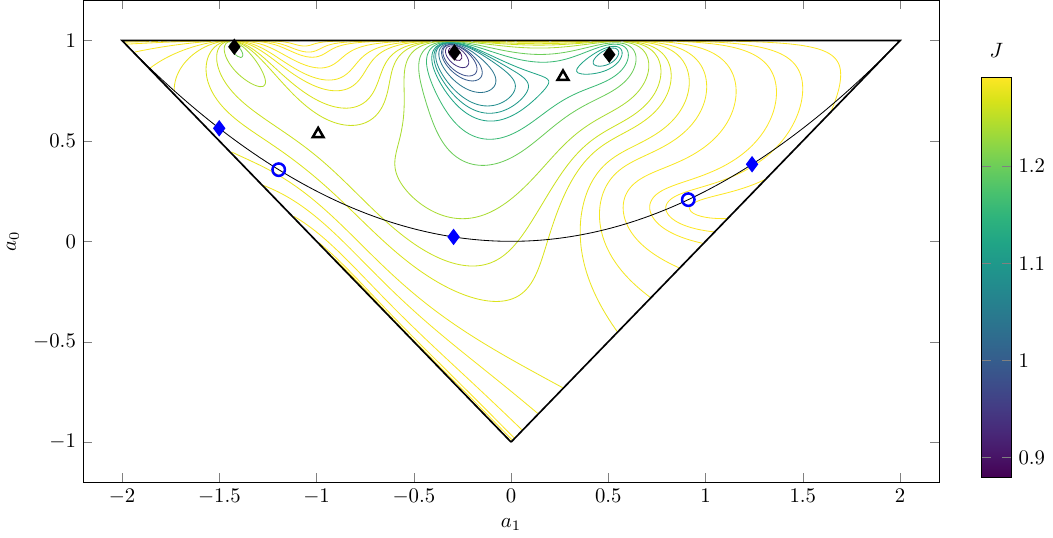}
    \caption[Contour plots of $\norm{H - \widehat H}_{\mathcal{H}_2}$ for varying reduced order models.]{
        Contour plot of the objective function value $\norm{H - \widehat H}_{\mathcal{H}_2}$ for varying reduced order models $\widehat{H}$.
        The problem setup is described in~\Cref{H2MOR_ex:didactic,exmp:enforcedMult}.
        Each point $(1:a_1:a_0)\in \Pset^2$ uses the optimal numerator $b$, obtained from solving~{\eqref{eq:optimal-b}}.
        The local minimizers (\textcolor{black}{$\blacklozenge$}) and saddle points (\textcolor{black}{$\triangle$}) in the displayed region are highlighted. 
        The black line depicts those points for which the discriminant $\Delta(a)$ vanishes. 
        When the model reduction problem is restricted to this discriminant, a different set of local minimizers (\textcolor{blue}{$\blacklozenge$}) and maximizers (\textcolor{blue}{$\bigcirc$}) is obtained. 
        }
    \label{intro_fig:mor}
\end{figure}

\begin{example}[Discrete-time]\label{H2MOR_ex:didactic}
   We consider the example given in \cite[Example~2]{Lagauw2023}, in which the discrete-time model 
   \begin{equation}
    H(z) = \frac{0.0448z^5 + 0.2368z^4 + 0.0013z^3 + 0.0211z^2 + 0.2250z + 0.0219}{z^6 -1.2024z^5 + 2.3675z^4 -2.0039z^3 + 2.2337z^2 -1.0420z + 0.8513}, \label{eq:discreteHigherorderModel}
   \end{equation}
   is approximated by a second-order model $\widehat{H}(z) = \frac{b_1z+b_0}{a_2z^2+a_1z+a_0}$.
   The authors report that the critical equations admit 49 isolated solutions, matching the prediction from \Cref{thm:numberSolutions}.
   For this problem instance, 11 critical points are real-valued, 5 of which correspond to stable models.

    For a fixed stable $a$ there is a unique $b$ that minimizes the error $J^2$.
    The derivation of the critical equations in \Cref{sec:CritEquations}, when carried out for $b$ only, shows that this optimal $b$ satisfies $ad - bc = \tilde{a} q$ for some $q \in \mathbb{C}[z]_{\le N-1}$, with $\tilde a(z) = z^na(1/z)$ as in \Cref{sec:discretetime}.
    The optimal $b$ can thus be found from the unique solution of the linear system
    \begin{equation}\label{eq:optimal-b}
      \begin{bmatrix} \Toeplitz{c}{n}\tran & \Toeplitz{\tilde a}{N}\tran \end{bmatrix}
      \begin{bmatrix} \vect{b} \\ \vect{q} \end{bmatrix} = \Toeplitz{a}{N}\tran \vect{d}.
    \end{equation}
    Its coefficient matrix, the Sylvester matrix of $c$ and $\tilde a$, is invertible
    for all points of interest because $c$ and $\tilde a$ have no common roots when $H$ and $\widehat H$ are stable.
    
    This observation allows us to draw a contour plot  of the objective function value $J^2$ for various reduced models $\widehat{H}$, shown in \Cref{intro_fig:mor}.
    We work in the affine chart $a_2=1$ and the plot over $(a_1,a_0) \in \RR^2$ is restricted to the triangular region corresponding to stable $\widehat{H}$; in the discrete-time setting, stability means that all roots of $a$ lie strictly inside the complex unit disk.\qedex
\end{example}

We now turn our attention to approximation problems where the bad locus $\mathcal{B}$ is not empty and hence the solution count has to be adapted. Such examples are mainly of theoretical interest and illustrate the assumptions and bounds obtained in the theorems of \Cref{sec:finitAndDeg,sec:otherConsid}.

\begin{example}[Non-empty bad locus]\label{ex:mulitpole_higher_mult}
    Consider the third-order model ($N=3$)
    \[
    H(s) = \frac{d}{c} = \frac{s^2 - 9 s + 20}{s^3 + 7s^2 + 16s + 12} = \frac{(s-4)(s-5)}{(s+2)^2(s+3)},
    \]
    which has two distinct poles, one of which has multiplicity two. 
    We compute its optimal first-order approximant ($n=1$). 
    The determinantal variety $\mathcal{D}$ contains the following approximate values of $[a] = (a_0:a_1) \in \Pset^1$:
    \[
    \mathcal{D} = \Set{ (-2.6300 : 1), \, (-2.0000 : 1), \, (0.3021 : 1), \, (4.4588 : 1), \, (33.8690 : 1)}.
    \]


    The generic count from \Cref{tab:numSols} lists $5$ solutions for $(N,n) = (3,1)$. 
    In our example, however, $c$ fails the genericity condition because it has a double root.
    Since $\mathcal{D}$ is finite, we may apply \Cref{H2MOR_TP:multiplerootscount} to adjust the generic count: we must subtract at least $\binom{N-2}{n-1} = \binom{1}{0} = 1$ to account for the non-empty bad locus $\mathcal{B}$.
    Indeed, symbolic computation shows that $\mathcal{B}$ consists of the single point $(-2:1)$.
    The adjusted count yields $4$ critical points. 

    The bad locus point can also be read off from \Cref{lem:MORbmatrix}.
    By the lemma, $[a] \in \mathcal{B}$ if and only if $\deg\gcd(c^\h, (\tilde a^\h)^2) > 1$. 
    Since $\deg c = 3$, $c$ has no roots at infinity, so $[a] \in \mathcal{B}$ requires $(s+2)^2 \mid \tilde a^2$, i.e., $\tilde a$ is proportional to $s + 2$ and hence $a = a_1(s - 2)$ for some $a_1 \in \Cset$. \qedex  
\end{example}

\begin{example}[Multiplicity in the bad locus points]\label{ex:badlocusmult}
    Reconsider the third-order ($N=3$) model $H$ from \Cref{ex:mulitpole_higher_mult}, but now compute its optimal second-order approximant ($n=2$).
    In this case, the determinantal variety $\mathcal{D}$ consists of 6 isolated points, one of which---the point $(a_0 : a_1 : a_2) = (6 : -5 : 1)$---occurs with multiplicity two with respect to the ideal of maximal minors of $\matx{M}(a,c,d)$.
    Hence, $\# \mathcal{D} = 6$ but $\deg \mathcal{D} = 7$, the latter of which coincides with the entry in \Cref{tab:numSols} for $(N,n) = (3,2)$.

    By \Cref{lem:MORbmatrix}, $[a] \in \mathcal{B}$ if and only if $\deg\gcd(c^\h, (\tilde a^\h)^2) \geq 3$. 
    Because $\deg c^\h = 3$, this occurs only when $c^\h$ divides $(\tilde a^\h)^2$, which requires $\tilde a$ to be proportional to $(s+2)(s+3)$, i.e., $a = a_2(s-2)(s-3) = a_2(s^2 - 5s + 6)$.
    Hence $\mathcal{B}$ consists of the single point $[a] = (1:-5:6)$.
    It has multiplicity two with respect to the maximal minors of $\matx{B}(a,c)$ as well, so that $\# \mathcal{B} = 1$ and $\deg \mathcal{B} = 2$.
    We conclude that there are $\# \mathcal{D} - \# \mathcal{B} = 5$ critical points.
    
    \Cref{H2MOR_TP:multiplerootscount} adjusts the generic count of $7$ by subtracting $\binom{N-2}{n-1} = \binom{1}{1} = 1$, yielding an upper bound of $6$ critical points. 
    The sharper bound from \Cref{cor:multbound}, which takes the multiplicity of the bad locus points into account, gives $\deg \mathcal{D} - \deg \mathcal{B} = 7 - 2 = 5$. \qedex
\end{example}

    \begin{example}[Critical point with multiplicity]\label{ex:critMult}
        Consider the higher-order model ($N=2$) defined by 
        \[
        H(s) = \frac{11s+13}{(s+7)(s+\frac{1}{5})}
        \]
        and consider reduction to order $n=1$. 
        Because $c$ has distinct roots and $\gcd(c,d) = 1$, the bad locus is empty and the determinantal variety $\mathcal D$ is finite (\Cref{lem:MORbmatrix} and \Cref{cor:finiteness_by_simple_roots}, respectively). 
        The matrix $\matx M(a,c,d) \in \CC[\vect a]^{3\times 3}$ is square and its determinant is of degree $3$. It has a double root at $a_0/a_1 = 1$ and a single root at $a_0/a_1 = \frac{91}{55}$. 
        Hence, there are $m = \#\mathcal{D} =2$ critical points, both real and stable. 
        Because of the root with multiplicity, the inequality in \eqref{eq:numSolutions} is strict: $2 = \#\mathcal{D} < \deg \mathcal{D} = 3$. 
        
        These input data $(c,d)$ are thus not generic in the sense of the fiber dimension theorem (\Cref{thm:fiber_dimension}): they lie outside the dense open subset on which the fibers of $\tau$, from \eqref{eq:tau}, are reduced.
        Hence the conclusion from \Cref{cor:finiteness} does not hold and the fiber is \emph{not} reduced.     
        \qedex
    \end{example}
    
Note that $(c,d)$ were deliberately chosen to have two critical points coincide. 
Such input data lie on a hypersurface in $CD$, which is a set of measure zero, and an arbitrarily small perturbation of $(c,d)$ splits the double critical point into two simple ones. 

\begin{example}[Prescribed multiplicity structure]\label{exmp:enforcedMult}
    Reconsider the setup from \Cref{H2MOR_ex:didactic}. 
    The black curve in \Cref{intro_fig:mor} is the zero-set of the discriminant $\Delta(a) = a_1^2 - 4 a_0 a_2$, depicting all reduced models with a multiple pole. 
    It is the submanifold described in~\Cref{rem:multiplicity}.
    
    The local minimizers of $J^2$, when restricted to this curve, satisfy (the discrete-time counterpart of) the conditions of~\Cref{thm:interpolation} (Meier--Luenberger), rather than those of \Cref{H2MOR_thm:walsh} (Walsh). 
    Indeed, when one parametrizes $\widehat H$ with a pole of multiplicity $m=2$, there remain only $m + 1 = 3$ interpolation conditions at $1/p$, instead of the $2m = 4$ conditions from \Cref{H2MOR_thm:walsh}. 
    Hence, when parameterizing $a \coloneq a_2(z-p)^2$, the factor $\tilde a^2 = (1-pz)^4$ in~\eqref{eq:lofs} is to be replaced by $\tilde a' \coloneq (1-pz)^3$:
    \begin{equation}\label{eq:MLcritical}
      a d - b c = (1-pz)^3 h, \qquad h \in \mathbb R[z]_{\le N-2}.
    \end{equation} 
    Upon normalization (e.g., $a_2 = 1$), these are $8$ polynomial equations in the $8$ unknowns $p, \vect{b}$ and $\vect{h}$.
    The system~\eqref{eq:MLcritical} is equivalent to
    \[
    \matx M'(p) \pmat{-1\\\vect b\\\vect h} = \vect{0}, \qquad \matx M'(p) \coloneq \bmat{\Toeplitz{a}{N}\tran \vect d &  \Toeplitz{c}{2}\tran &  \Toeplitz{\tilde a'}{N-1}\tran } \in \mathbb R[p]^{8\times8}.
    \]
    This matrix is square and depends on the single parameter $p$; it is a polynomial eigenvalue problem (cf.~Example~\ref{ex:PEP}). 

    For the higher-order model in \eqref{eq:discreteHigherorderModel}, $\det \matx M'(p)$ has $3N - 1 = 17$ distinct roots, $7$ of which are real, and $3$ of which lie in the unit disk: $p \in \{-0.6196, 0.1478, 0.7503\}$.\footnote{Similar to $\matx M(a,c,d)$, the matrix $\matx M'(p)$ may have a bad locus, consisting of the values of $p$ for which its last $N+1$ columns, $[\Toeplitz{c}{2}\tran \, \Toeplitz{\tilde a'}{N-1}\tran]$, drop rank. 
    Its behavior can be studied similarly as in \Cref{lem:MORbmatrix}.}
    These constitute the three local minimizers of the restricted problem, depicted in blue in \Cref{intro_fig:mor}. 
    
    The restricted problem has two local maxima, $p \in \{-0.4556,0.5974\}$, which do not appear as a solution to \eqref{eq:MLcritical}. 
    At both points, the optimal numerator $b$ is such that $a$ and $b$ share a common factor---the stationary point is \emph{degenerate}---so that $\widehat H = b_1/(z-p)$ has McMillan degree one instead of two.
    These local maximizers are stationary points of the (unrestricted) first-order problem ($n=1$).    \qedex
\end{example}

\begin{remark}[Degenerate stationary points]\label{rem:degen}
    Also in the unrestricted case, stationary points of a lower-order problem can reappear in the higher-order problem as degenerate stationary points, see e.g.,~\cite[Proposition~2.7]{Baratchart1992} for a characterization in discrete time. 
    Because such points lie outside of $\mathcal{M}_n$, \Cref{H2MOR_thm:walsh} does not apply to them, and they need not appear as a solution of the critical equations \eqref{eq:criticalEqsMatrix}.
    For model reduction purposes, this poses no problem as~\cite[Theorem~5.3]{Regalia1995} guarantees that such degenerate points can never be local minimizers.
\end{remark}

\begin{example}[Optimal reduced model with a triple pole]\label{exmp:vandooren}
    We revisit the example of \cite[Section~6]{VanDoren2010}, following their Theorem~6.2, where
    \[
      H(s) = \frac{\frac{1}{256}(s^3 - 9s^2 + 295s + 2449)}{s^4 + 13s^3 + \frac{513}{16}s^2 + \frac{239}{8}s + \frac{153}{16}}
    \]
    has four distinct poles. 
    For $n = 3$, the bad locus $\mathcal B$ is therefore empty, and the determinantal variety $\mathcal D$ consists of $15 = 2^4 - 1$ points of multiplicity one, in agreement with~\Cref{thm:numberSolutions}. 
    Three of these points are real, and only one of them corresponds to a stable model: $a = (s+1)^3$ with $b = 1$. 
    Since a global minimizer exists (cf.~\Cref{rem:existence}), $\widehat H = 1/(s+1)^3$ is the global minimizer.
    
    For this particular higher-order model $H$, the global minimizer thus lies on the curve of polynomials $a = (s-p)^3$ with a triple root, which is contained in the discriminant locus $\Delta(a) = 0$. 
    Being a stationary point of the unrestricted problem, $\widehat H$ satisfies~\Cref{H2MOR_thm:walsh}: the interpolation conditions~\eqref{eq:interpolationConditions_Cont} hold at $s = 1$ for $k = 0, \dots, 5$. 
    These include the conditions for $k = 0, \dots, 3$ required by~\Cref{thm:interpolation}, so that the latter are automatically satisfied as well.
    
    Hence, this stationary point would also have been obtained by enforcing the multiplicity structure from the start, as in~\Cref{exmp:enforcedMult}, i.e., by replacing $\tilde a^2 = (s+p)^6$ in~\eqref{eq:lofs} with $(s+p)^4$. 
    One may verify that in this case, $\det \matx M'(p)$ has $15$ distinct roots, $5$ of which are real and $3$ of which are stable: $p \in \{-1, -0.084, -0.476\}$. 
    Contrary to the point $p = -1$, which corresponds to the global minimizer of the unrestricted problem, the other two are stationary points of the restricted problem only.\qedex
\end{example}

\section{Conclusions and outlook} \label{sec:Conclusions}

We have illustrated how tools from algebraic geometry can facilitate the analysis of systems of polynomial equations arising in engineering applications. 
In the specific context of the $\Hardy{}$-norm model order reduction problem, we applied the fiber dimension theorem and the Thom--Porteous formula to obtain an exact a priori count of the complex solutions to the critical equations. 
This result is substantially sharper than the coarse bounds available from, for instance, Bézout's theorem. 
Our formula shows that the number of critical points grows rapidly with the problem size, reflecting the computational complexity observed in engineering practice.

The most demanding part of the analysis was proving that the incidence variety $\mathcal{W}$ defined by the critical equations is irreducible, reduced and of the same dimension as the space of input variables.
By contrast, the analysis of the bad locus---corresponding to extraneous solutions introduced by the MEP formulation---was favorable, as we were able to show that the bad locus is most often empty. 

The methodology presented here is not restricted to the model reduction problem under consideration, but can be applied more broadly to other polynomial systems and multiparameter eigenvalue problems arising in related applications. 
For instance, in our related work~\cite{Kayser2026a} we employed a similar strategy in the context of a system identification problem.

In \Cref{ex:mulitpole_higher_mult,ex:badlocusmult}, we illustrated that for specific problem instances, in particular when the higher-order model has multiple roots, the bad locus $\mathcal{B}$ can be non-empty and contain points of higher multiplicity.
Moreover, \Cref{exmp:enforcedMult} touches on how the critical equations change when the poles of the reduced model follow a prescribed multiplicity structure. 
A detailed understanding of these special cases is an interesting direction for future work. 

In the case of a higher order model of full McMillan degree $N$ and simple poles, we have reproven the main result of \cite{Hanzon2012} on finiteness of critical points (cf.~\Cref{cor:finiteness_by_simple_roots}). It remains an open question whether the critical points can ever form an infinite set, for \emph{any} higher order model (of degree $N\geq n$). We conjecture that the critical points \emph{always} form a finite set, and have verified this in the case of a finite bad locus.

\Cref{exmp:enforcedMult} also illustrated that stationary points of a lower-order problem can reappear in a (possibly restricted) higher-order problem as \emph{degenerate} stationary points, that is, as models with pole-zero cancellations and thus McMillan degree strictly smaller than $n$. 
For the (unrestricted) problem in discrete time, such points were characterized in~\cite[Proposition~2.7]{Baratchart1992}. 
By formulating the model order reduction problem \eqref{eq:modelReductionOptProblem} over $\mathcal{M}_n$, we have excluded these degenerate points from the search space. 
However, even though they are never local minimizers~\cite[Theorem~5.3]{Regalia1995}, they can appear as saddle points or local maxima in the optimization landscape and may therefore still be of practical interest. 
A systematic understanding of how many such degenerate stationary points arise, and how they influence the optimization landscape, is another interesting research direction. In our related work~\cite{Kayser2026a}, 
lower-order solutions similarly appear as separate components of the critical equations.

Another natural extension is to investigate how the presented analysis generalizes to frequency weighted $\Hardy$-norm model order reduction (see, e.g.,~\cite{Huang2001}), where a nontrivial weight function is incorporated in the approximation criterion.
We expect the weights to alter the structure of the critical equations, and thereby to complicate the analysis of the Walsh variety $\mathcal{W}$ and corresponding bad locus.

\subsection*{Acknowledgements}

We would like to thank Simon Telen and Fulvio Gesmundo for helpful discussions and feedback on this manuscript. 
Parts of this article are based on the authors' PhD theses: \cite[Chapter 5]{Leonie2025} and \cite[Chapter 5]{Sibren2025}. 

This work was supported in part by the KU Leuven Research Fund (projects iBOF/24/064 and IDN/24/013), the Industrial Research Fund and several LRD bilateral industrial projects; in part by the Flemish Government agencies: FWO (TBM T001725N), EWI through the Flanders AI Research Program and VLAIO (CSBO - HBC.2021.0076); and in part by the European Commission (ERC Adv. Grant) under Grant 885682. Views and opinions expressed are however those of the author(s) only and do not necessarily reflect those of the European Union or ERC. Neither the European Union nor the granting authority can be held responsible for them. The work of S.~Lagauw was supported by the FWO Fundamental Research Fellowship under Grant 11K5625N.

\subsection*{Declaration of AI use}
GenAI-based tools were used solely in the proofreading process.

\appendix

\section{Relevant background in Algebraic Geometry}\label{app:AG}

This appendix reviews the algebro-geometric concepts that lie at the basis of our analysis. 
The exposition is intentionally limited in scope and presents only the essential ideas for our approach.
For readers interested in applications, \cite{Sommese2005} provides an excellent introduction to polynomial system solving in the sciences and engineering. The book \cite{Clader2025} also provides a beginner-friendly introduction to algebraic geometry.
Basic results in algebraic geometry are introduced in e.g., \cite{Harris1995,Shafarevich2013,Cox2015}, while more advanced treatments can be found in \cite{Hartshorne1977,EisenbudHarris2016}. 

\subsection{Affine varieties}\label{AG:varieties}

Let $R_n \coloneq \Cset[x_1,\dots,x_n]$ denote the ring of (multivariate) polynomials in the variables $x_1,\dots,x_n$.
An \emph{affine variety} $X \subseteq \Cset^n$ is the set of all points in $\Cset^n$ at which a given collection of polynomials $F = \{f_1, \dots, f_s\} \subseteq R_n$ vanishes simultaneously
\[
X = V(F) \coloneq\Set{\vect{x} \in \Cset^n | f(\vect{x}) = 0 \text{ for all }f \in F}, \qquad F \subseteq R_n. 
\]
\emph{Closed subsets} of $\Cset^n$ are precisely the sets of the form $V(F)$ for some $F \subseteq \Cset[x_1,\dots,x_n]$; their complements are called \emph{open sets}.

A \emph{subvariety} or \emph{closed subset} $Y$ of $X$ is an affine variety $Y \subseteq X$; an \emph{open subset} is of the form $X \setminus Y$. 
The closure $\overline{Y}$ of a subset $Y \subseteq X$ is the smallest closed subset $\overline Y \subseteq X$ containing $Y$. 
The set $Y$ is \emph{dense} in $X$ if $\overline Y = X$.
A property holds for \emph{generic} points $\vect{x} \in X$ if it holds for all points in a dense open subset of $X$, so that \enquote{almost all} points of $X$ share that property.
The variety $X$ is \emph{irreducible} if $X \neq \emptyset$ and whenever $X = Y_1\cup Y_2$ with subvarieties $Y_1,Y_2 \subseteq \Cset^n$, then $X=Y_1$ or $X=Y_2$. 
Equivalently, every two non-empty open subsets of $X$ intersect.

These closed sets, as defined above, endow $\Cset^n$ with the so-called \emph{Zariski topology} (see, e.g., \cite{Shafarevich2013} for a formal treatment).
This topology is much coarser than the usual (Euclidean) topology on $\Cset^n$; for instance, any non-empty Zariski-open subset of an irreducible variety $X$, e.g., $X = \Cset^n$, is dense. 
An open subset of an affine variety is a \emph{quasi-affine variety}; almost all of the definitions in this section apply also to such varieties. 

Every variety $X$ can be uniquely written as a finite union of its \emph{irreducible components}, in such a way that no component is contained in the union of the others (an \emph{irredundant} decomposition):
\[
X = X_1 \cup\dots \cup X_r, \qquad X_i\text{ irreducible, closed subset }X_i\not\subseteq X_j \text{ for all }i \neq j.
\]

The \emph{dimension} of an irreducible variety $X$ is the length of the longest strictly increasing chain of irreducible subvarieties
\[
\dim X = \max \Set{d \in \Nset | \exists Y_0 \subsetneq Y_1 \subsetneq  \dots \subsetneq Y_d = X, \text{ $Y_i$ irreducible subvarieties}}.
\]
By convention $\dim \emptyset \coloneq -1$.
Let $X$ be an arbitrary variety, and let $X = X_1 \cup \cdots \cup X_r$ be its decomposition into irreducible components. Then,
\[
\dim X \coloneq \max_{1 \le i \le r} \dim X_i.
\]
It turns out that any chain of irreducible varieties can be extended (by suitably inserting more irreducible subvarieties) to a maximal one, and if $X$ is irreducible, then all such maximal chains have the same length. 
In particular, any closed proper subvariety of an irreducible variety has strictly smaller dimension than the variety itself.

A variety $X$ is \emph{equidimensional} if every irreducible component $X_i \subseteq X$ has the same dimension.
The \emph{codimension} of an irreducible subvariety $Y$ in a variety $X$ is $\codim_X(Y) = \dim X' - \dim Y$, where $X' \subseteq X$ is an irreducible component of largest dimension containing $Y$. 
The dimension at a point $\vect{p} \in X$ is $\codim_X\{\vect{p}\}$; if $X$ is irreducible then this is just the dimension of $X$.

\subsection{Reducedness and smoothness}\label{sec:smoothness}
The \emph{defining ideal} of an affine variety $X \subseteq \CC^n$ is the set of all polynomials vanishing on $X$,
\[
I(X) = \Set{f \in R_n | f(\vect{x}) = 0 \text{ for all } \vect{x} \in X}.
\]
For a set of polynomials $F \subseteq R_n$, denote by $\langle F \rangle$ the ideal generated by $F$.
Conversely, for an ideal $I \subseteq R_n$, a subset $F = \{f_1,\dots,f_s\} \subseteq I$ is called a \emph{generating set} of $I$ if $I = \langle F \rangle$.

For any ideal $J \subseteq R_n$, its \emph{radical} $\sqrt{J}$ is defined by
\[
\sqrt{J} \coloneq \Set{f \in R_n | \exists \, k \in \Nset: f^k \in J}.
\]
A set of polynomials $F \subseteq R_n$ is said to define a variety $X$ \emph{in a reduced way} if $V(F) = X$ and the ideal $\langle F \rangle$ is radical, i.e., $\sqrt{\langle F \rangle} = \langle F \rangle$. 
By \emph{Hilbert's Nullstellensatz} we have
\[
I(V(F)) = \sqrt{\langle F \rangle},
\]
so if $V(F) = X$, then $F$ defines $X$ in a reduced way if and only if $\langle F \rangle = I(X)$, that is, if $F$ is a \emph{generating set} of the defining ideal $I(X)$.

If $F = \{f_1,\dots,f_m\}$ is a generating set of the defining ideal $I(X)$, then the \emph{smooth locus} of $X$ is
\[
X^{\text{sm}} = \Set{ \vect{p} \in X \mid \rank \matx{J}_F (\vect{p}) \geq n - \codim_X \{\vect{p}\}},
\]
where $\matx{J}_F(\vect{p})$ is the Jacobian matrix of $F$ at $\vect{p}$:
\[
\matx{J}_F(\vect{p}) = \left[\frac{\partial f_i}{\partial x_j}(\vect{p})\right]_{i,j}.
\]
The smooth locus is an open dense set $X^{\text{sm}} \subseteq X$, and, if $X$ is irreducible, is a connected complex manifold of the same (complex) dimension as $X$. 
A point $\vect{p} \in X$ is called \emph{singular} if it does not lie in $X^{\mathrm{sm}}$.
We denote $X^{\text{sm}}(\Rset) \coloneq X^{\text{sm}} \cap \Rset^n$, that is, the set of smooth points of $X$ with real coordinates.
If the defining polynomials are real and $\vect{p} \in X^{\text{sm}}(\Rset)$, then $X^{\text{sm}}(\Rset) \subseteq \Rset^n$ is a smooth (real) submanifold in a neighborhood of $\vect{p}$.

For polynomial equations $F$ we defined what it means to defined $X=V(F)$ in a reduced way. 
We now turn from this global viewpoint to a local one. 
Globally, we consider polynomial functions defined on all of $\CC^n$, and reducedness is expressed by the defining ideal being radical.
Locally, however, we focus on the behavior near a single point. 
Functions in this local setting are defined only on a small neighborhood and are described algebraically by localizations of rings, which allow denominators that do not vanish at the point. 
To formulate the local analogue of reducedness, we therefore begin by clarifying what we mean by such  local functions.

Let $X \subseteq \mathbb{C}^n$ be an affine variety and $\vect p \in X$. A regular function germ $\varphi_{\vect p}$ at $\vect{p}$ is a function $\varphi \colon U \to \CC$ defined on some open set $\vect p \in U \subseteq X$ which, on $U$, can be represented as $\varphi = \frac{f}{g}$, $f,g \in R_n$ and $V(g) \cap U= \emptyset$. Two such function germs are considered equal if they agree on a common open neighborhood of $\vect p$.
That is, germs are functions on an arbitrarily small neighborhood of $\vect p$, so it does not make sense to evaluate $\varphi_{\vect p}$ at any point of $X$ other than $\vect p$, though local properties like order of vanishing, gradient and such are still well-defined. The collection of regular function germs form the \emph{local ring of $X$ at $\vect p$} denoted by $\mathcal{O}_{X,\vect p}$. For example, regular function germs on $X = \CC^n$ at $\vect 0$ are rational functions without a pole in $\vect 0$.

Every polynomial $f \in R_n$ defines function germs $f_{\vect p} \in \mathcal{O}_{\CC^n,\vect p}$. The \emph{germ} of an ideal $\langle F \rangle_{R_n}$ is the ideal $F\mathcal{O}_{\CC^n,\vect p}$ in the ring $\mathcal{O}_{\CC^n,\vect p}$ generated by $\set{f_{\vect p} | f \in F}$. Now $F$ defines the variety $X = V(F) \subseteq \CC^n$ \emph{in a reduced way at $\vect{p}$} if $F\mathcal{O}_{\CC^n,\vect p} = I(X)\mathcal{O}_{\CC^n,\vect p}$, equivalently if $F\mathcal{O}_{\CC^n,\vect p} $ is a radical ideal.
For example, if the Jacobian matrix of $F$ at $\vect{p}$ is of maximal rank, then $F$ defines $X=V(F)$ in a reduced way at $\vect{p}$ and furthermore $\vect{p}$ is a smooth point of $X$.

\subsection{Dimension of fibers and intersections}\label{sec:dimFibers}
 
A \emph{morphism} of affine varieties $X \subseteq \Cset^n$ and $Y \subseteq \Cset^m$ is a coordinate-wise polynomial map \[
\phi = (\phi_1,\dots,\phi_m)\colon \Cset^n\to \Cset^m, \qquad \phi(X) \subseteq Y.
\]
Polynomial coordinate functions is what sets a morphism apart from general maps; it ensures that the preimages of subvarieties are again subvarieties, and that the composition of $\phi$ with a regular function germ $\varphi_{\vect q}$, $\phi(\vect p) = \vect q$, yields another regular function germ $(\varphi \circ \phi)_{\vect p}$.
The morphism is \emph{dominant} if $\overline{\phi(X)} = Y$. 
The \emph{fiber} $\phi^{-1}(\vect{q})$ of $\vect{q} \in Y$ is the subset $\set{\vect{x} \in X | \phi(\vect{x}) = \vect{q}}$. 
Equations of $\phi^{-1}(\vect{q})$ are given by 
\[
I(X) \cup \{\phi_1-q_1,\dots,\phi_m-q_m\} \subseteq \Cset[x_1,\dots,x_n],
\]
and $\phi$ has \emph{reduced fibers at $\vect{q}$} if this set of equations defines $\phi^{-1}(\vect{q})$ in a reduced way.
An \emph{isomorphism} is a morphism $\phi \colon X \to Y$ that is bijective and whose inverse $\phi^{-1}\colon Y \to X$ is also a morphism; in this case $X$ and $Y$ are said to be \emph{isomorphic}.

The following two fundamental results of algebraic geometry can be seen as extensions of the rank-nullity theorem from linear algebra to non-linear equations.

\begin{theorem}[Fiber dimension theorem {\cite[Theorem I.1.25]{Shafarevich2013}, \cite[Corollary III.10.7]{Hartshorne1977}}]\label{thm:fiber_dimension}
Let $\phi \colon X \to Y$ be a morphism of varieties and $X$ irreducible.
\begin{enumerate}
\item For every $\vect{q} \in \phi(X)$, one has $\dim \phi^{-1}(\vect{q}) \geq \dim X - \dim \overline{\phi(X)}$.
\item There is a dense open subset $U \subseteq \overline{\phi(X)}$ such that $\phi$ has reduced (even smooth!) fibers at $\vect{q} \in U$ and equality holds in the first statement.
\end{enumerate}
\end{theorem}

\begin{theorem}[Dimension of intersection {\cite[Theorem I.1.24]{Shafarevich2013}}]\label{thm:intersection_dim}
\leavevmode
\begin{enumerate}
\item If $X,Y\subseteq \Cset^n$ are equidimensional varieties, then every irreducible component $Z$ of the intersection $X \cap Y$ satisfies
\[
\dim Z \geq \dim X + \dim Y - n.
\]
\item If $f_1,\dots,f_s \in R_n$, then $V(f_1,\dots,f_s)$ is empty or each of its irreducible components has dimension $\geq n-s$.
\end{enumerate}
\end{theorem}

With minor additional technicalities (defining morphisms of quasi-affine varieties), the theorem remains valid after replacing $\Cset^n$ by any non-empty open subset $U \subseteq \Cset^n$.

Consider a non-empty open subset $U \subseteq \Cset^n$ and a variety $X \subseteq U$ with $\codim_U X = s \leq n$. 
We say that $X$ is a \emph{set-theoretic complete intersection in $U$} if there exists a set $\{f_1,\dots,f_s\} \subseteq R_n$ such that $X = V(f_1,\dots,f_s) \cap U$. 
In this case, $X$  can be defined by $s$ equations, and no strictly smaller number of equations suffices to cut out $X$ set-theoretically (by \Cref{thm:intersection_dim}).
If, in addition, $f_1,\dots,f_s$ generate the defining ideal $I(X)$ of $X$ in $\mathcal{O}_{\CC^n,\vect p}$ (i.e., they define $X$ by a radical ideal in $\mathcal{O}_{\CC^n,\vect p}$) for all $\vect p \in U$, then $X$ is an \emph{ideal-theoretic complete intersection in $U$}.

\begin{theorem}[Unmixedness theorem {\cite[Corollary 18.14]{Eisenbud1995}}]\label{thm:unmixedness}
Let $X = V(f_1,\dots,f_s)$ be a set-theoretic complete intersection.
If $X$ is defined in a reduced way somewhere on every irreducible component of $X$, then it is reduced everywhere, so actually an ideal-theoretic complete intersection.
\end{theorem}

We now introduce affine bundles as a special and useful class of subvarieties.
Take a (quasi) affine variety $B \subseteq \Cset^m$.
Fix an $(\ell \times n)$-matrix $\matx{M}(\vect{y})$ whose entries are polynomials in the coordinates $\vect{y} \in B$. 
Assume that $\rank \matx{M}(\vect{y}) = k$ at every point $\vect{y} \in B$. 
Choose a polynomial column vector $\vect{b}(\vect{y}) \in \Cset[\vect{y}]^{\ell}$ so that for all $\vect{y} \in B$, the vector $\vect{b}(\vect{y})$ lies in the image of $\matx{M}(\vect{y})$, ensuring that the linear system $\matx{M}(\vect{y})\vect{x} = \vect{b}(\vect{y})$ admits solutions for every $\vect{y} \in B$.
Under these assumptions, let
\[
X = \Set{ (\vect{x},\vect{y}) \in \Cset^n \times B |\matx{M}(\vect{y})\vect{x} = \vect{b}(\vect{y})}
\]
and consider the projection $\pi: X \to B, (\vect{x},\vect{y}) \to \vect{y}$.
For a chosen point $\vect{y} \in B$, the fiber
\[
\pi^{-1}(\vect{y}) = \Set{ \vect{x} \in \Cset^n  |  \matx{M}(\vect{y})\vect{x} = \vect{b}(\vect{y})} \times \{\vect y\},
\]
is an affine linear subspace of $\Cset^n$ of dimension $r = n-k$, hence isomorphic to $\Cset^{r}$.
We say that $\pi: X \to B$ is an \emph{affine bundle of rank $r$}.
In that context, $X$ is called the \emph{total space} of the bundle, $B$ is the \emph{base space}, and $\pi$ is the \emph{bundle projection}.

Establishing that a given morphism $\pi: X \to B$ is a rank $r$ affine bundle is advantageous.
Because all fibers have dimension $r$, we have $\dim X = \dim B+r$.
Moreover, the local product structure of the bundle---namely, that $X$ is locally isomorphic to $B \times \Cset^r$---ensures that properties such as irreducibility, equidimensionality, reducedness, and smoothness pass from the base $B$ to the total space $X$. 
For instance, if $B$ is irreducible, then so is $X$. 

\subsection{Projective space and determinantal varieties}\label{AG:projective}

Projective space $\Pset^n$ is defined as $(\Cset^{n+1} \setminus \{0\})$ modulo the equivalence relation that identifies two nonzero vectors $\vect{x}, \vect{x}' \in \Cset^{n+1}$ whenever they are scalar multiples of each other.
A projective point $(x_0 : x_1 : \dots : x_n) \in \Pset^n$ represents the equivalence class $[\vect{x}]$ of all nonzero vectors $(\lambda x_0, \lambda x_1, \dots, \lambda x_n) \in \Cset^{n+1}$ with $\lambda \in \Cset \setminus \{0\}$.
The notation $\Pset(V)$ denotes the \emph{projectivization} of a vector space $V$, i.e., the set of nonzero vectors modulo scaling, in other words one-dimensional linear subspaces of $V$. 
In particular, $\Pset(\Cset^{n+1}) = \Pset^n$.

A \emph{projective variety} $X$ is the zero-set $V(F) \subseteq \Pset^n$ of a set of homogeneous polynomials $F = \{f_1,\dots,f_s\} \subseteq R_{n+1} = \CC[x_0,\dots,x_n]$. Notions such as open/closed subsets are defined analogously.
Its defining ideal $I(X)$ is defined as the defining ideal of the \emph{affine cone} $\widehat X = \Set{\vect{x} \in \Cset^{n+1} | [\vect{x}] \in X} \cup \{0\}$, and it is generated by all homogeneous polynomials in $R_{n+1}$ vanishing on $X$. 
If $X \subseteq \Pset^n$ is a projective variety and $\widehat{X} \subseteq \Cset^{n+1}$ denotes its affine cone, then $\dim \widehat{X} = \dim X + 1$ and hence the codimension is preserved:
\[
\operatorname{codim}_{\Pset^n}(X)
= n - \dim X
= (n + 1) - (\dim X + 1)
= \operatorname{codim}_{\Cset^{n+1}}(\widehat{X}).
\]

Projective space $\Pset^n$ can be decomposed as the disjoint union of the hyperplane at infinity $H_\infty \coloneq V(x_0)$ and its complement 
\[
U_0 = \{x_0 \neq 0\} = \Set{(1:x_1:\dots :x_n) | \vect{x} \in \Cset^n},
\]
which is called the affine part on $\Pset^n$. 
A \emph{standard affine chart} of a projective variety $X$ is an open subset of the form $X \setminus V(x_i)$.
Affine charts can be naturally identified with affine varieties. For example, in the case $X=\PP^n$, normalizing $x_0 = 1$ gives the isomorphism $U_0 \cong \Cset^n$.

In this way, we can associate to a projective variety $X \subseteq \Pset^n$ not contained in $H_\infty$ the affine variety
\[
X' = X \cap U_0 = \Set{\vect{x} \in \Cset^n | (1:x_1:\dots:x_n) \in X}, \qquad I(X') = I(X)|_{x_0=1}.
\]
Conversely, if $Y \subseteq \Cset^n$, then one can consider the closure $\overline{Y} \subseteq \Pset^n$. 
Its defining ideal is given by
\[
I(\overline{Y}) = \langle f^{\mathrm{h}}=\sum_{i=0}^{d}x_{0}^{\,d-i}f^{(i)} \, | \, f=f^{(0)}+\dots+f^{(d)} \in I(Y),\ \deg f = d \rangle, 
\] where in this case, $f^{(i)}$ denotes the component of $f$ of total degree $i$.
Note, however, that if $Y = V(f_1,\dots,f_s)$, then it is usually \emph{not} the case that $\overline{Y} = V(f_1^{\text{h}},\dots,f_s^{\text{h}})$, even when the $f_i$ define $Y$ in a reduced way; the latter projective variety may contain extra components in $H_\infty$.

A \emph{quasi-projective variety} is an open subset of a projective variety. 
Quasi-projective varieties include all previously introduced types of varieties. 
Indeed, if $W \subseteq \Cset^n$ is an affine variety, and $\overline{W} \subseteq \Pset^n$ denotes its projective closure, then $W = \overline{W} \cap \Cset^n$, so an affine variety $W$ is an open subset of the projective variety $\overline{W}$ and hence is quasi-projective. 
It is the complement of the hyperplane at infinity.

Let $X \subseteq \mathbb P^n$ be a quasi-projective variety. 
A function $\varphi \colon X \to \Cset$ is \emph{regular} if on every standard affine chart $X \setminus V(x_i)$ it defines a regular function germ at every point $\vect p \in X_i$ in the sense of \Cref{sec:smoothness}. Locally, though not necessarily globally, it can be written as a fraction of two homogeneous polynomials $\varphi = f(x)/g(x)$ of the same degree (which hence is well-defined for projective points).

Affine varieties carry many regular functions, for instance all polynomials. Projective varieties, by contrast, carry very few. 
The following result is an algebraic counterpart of the fact that a holomorphic function on a compact, connected complex manifold is constant (Liouville's theorem/the maximum modulus principle).

\begin{theorem}[{\cite[Corollary I.5.1.1]{Shafarevich2013}}]\label{thm:constant}
Every regular function on an irreducible projective variety is constant.
\end{theorem}

An important advantage of working in projective space is that some equations admit solutions that do not exist in affine space. 
For example, parallel lines in the affine plane intersect at $H_\infty$. 
Formally, this principle is captured by the following theorem.

\begin{theorem}[{\cite[Theorem I.1.24]{Shafarevich2013}}]\label{thm:projective_varieties_meet}
If $X,Y \subseteq \Pset^n$ are projective varieties with $\dim X + \dim Y \geq n$, then $X \cap Y \neq \emptyset$ and $\dim X \cap Y \geq \dim X + \dim Y - n$.
\end{theorem}

\emph{Determinantal varieties} are varieties defined by the vanishing of the minors of a matrix.
In this work, we focus on determinantal varieties associated with matrices $\matx{M}(\vect{x}) \in R_{n+1}^{f\times e}$  ($e \leq f$) of homogeneous polynomials in $\vect{x} \in \Pset^n$, whose columns have constant degree $d_1,\dots,d_e$:
\begin{equation}
\matx{M}(\vect{x}) = (h_{i,j})_{1 \leq i \leq f, \; 1 \leq j \leq e}, \qquad h_{i,j} \in R_{n+1}  \text{ is homogeneous of degree } d_j. \label{eq:M}
\end{equation}
The determinantal variety of rank $k$ associated to $\matx{M}$ is the projective variety
\[
D_k(\matx{M}) = \Set{\vect{x} \in \Pset^n | \rank \matx{M}(\vect{x}) \leq k}.
\]
It is a closed set for which the equations are given by the $(k+1)$-minors of $\matx{M}$ which are homogeneous by assumption on the degrees of the entries of $\matx{M}$.

\begin{theorem}[Dimension of determinantal varieties {\cite[Theorem 14.4(b)]{Fulton1998}}]\label{thm:detlowerbound}
Let $\matx{M}$ be the matrix polynomial from \eqref{eq:M}. 
Every irreducible component $Z \subseteq D_k(\matx{M})$ one has
\[
\dim Z \geq n - (e-k)(f-k).
\]
\end{theorem}

One can show, if the expected dimension $n - (e-k)(f-k)$ is non-negative and under suitable positivity assumptions on the degrees $d_j$, that $D_k(\matx M)$ is always non-empty, in analogy to \Cref{thm:projective_varieties_meet}.
Note that for $e = 1$ and $k = 0$, one obtains a projective analogue of \Cref{thm:intersection_dim}(ii).

A helpful way to interpret \Cref{thm:detlowerbound} is to compare its conclusion with the \emph{expected} geometry of determinantal varieties. 
In the space of all $(f \times e)$-matrices, the locus of matrices of rank at most $k$ is defined by the $(k+1)$-minors and has codimension $(f-k)(e-k)$, hence dimension $k(f+e-k)$. 
One can see this dimension directly by parametrizing a rank-$k$ matrix via a $k$-dimensional row-space, a $k$-dimensional column space, and an invertible $(k\times k)$-core, which together give $k(f+e-k)$ degrees of freedom. 
When the entries of the matrix vary independently, this ``expected'' dimension is achieved exactly. 
In contrast, when the entries arise from prescribed homogeneous polynomials (e.g., the matrix $\matx{M}(\vect{x})$ from \eqref{eq:M}) the rank conditions (potentially) impose fewer independent constraints, and the codimension can only become smaller.
\Cref{thm:detlowerbound} formalizes this: every component of $D_k(\matx{M})$ has codimension no greater than in the unrestricted case.

\subsection{Degree and a special case of Thom--Porteous}\label{AG:degreeTP}

If $X = V(F) \subseteq \Cset^n$ is a finite set (potentially empty), then the quotient space $R_n/\langle F \rangle$ has finite dimension. 
The degree of the ideal, $\deg \langle F \rangle$, is the dimension of this quotient space and we define the \emph{degree} of $X$ as $\deg X \coloneq \#X = \deg I(X)$.
In this case, $\deg \langle F \rangle$ corresponds to the number of points in $X$, counted with multiplicity, and we have
\[ 
\deg \langle F \rangle = \#X \quad \iff \quad \text{$F$ defines $X$ in a reduced way}.
\]
When $X$ is infinite (its dimension is positive) the definition above is no longer meaningful. 
In this case, we pick $\dim X$ generic affine-linear polynomials $\ell_1,\dots,\ell_{\dim X}$ and define 
\[
\deg X := \deg \langle I(X),\ell_1,\dots,\ell_{\dim X}\rangle = \#(X \cap V(\ell_1,\dots,\ell_{\dim X})).
\]
The degree is independent of the particular choice of the $\ell_i$, and this construction agrees with the familiar notion of the degree of a hypersurface: if $X=V(f)$ with $f$ square-free (e.g., irreducible), then $\deg X = \deg f$.

If no irreducible component of $X \subseteq \Pset^n$ is contained in the hyperplane at infinity $H_\infty$, then its degree is the degree of $X \cap \Cset^n$.
The degree (and dimension) of a projective variety can alternatively be characterized via the leading coefficient of the \emph{Hilbert function}:
\[
h_X(t) \coloneq \dim_\Cset \frac{S_t}{I(X)_t} = \frac{\deg X}{(\dim X)!}t^{\dim X} + O(t^{\dim X-1}),
\] 
where $S_t$ denotes the homogeneous polynomials in $R_{n+1}$ of degree $t$ and $I(X)_t$ the subspace of $S_t$ containing all elements vanishing on $X$. 
The term $O(t^{\dim X-1})$ is in the standard big-$O$ notation, meaning that the remainder is bounded by a constant multiple of $t^{\dim X-1}$ for all sufficiently large $t$.
This construction also gives a general definition of the degree of a homogeneous ideal $I \subseteq R_{n+1}$.
The degree of a non-empty projective variety is strictly positive, while the empty set has degree $0$ (its Hilbert function is the constant $0$).

A version of Bézout's theorem is that if $X = V(f_1,\dots,f_n) \subseteq \PP^n$ is finite, then
\[
\# X \leq \deg \langle f_1,\dots,f_n\rangle = \deg(f_1) \dotsm \deg(f_n).
\]
This has a vast generalization, of which we consider a special case. 
Recall from \Cref{thm:detlowerbound} that determinantal varieties have a general lower bound on their dimension. The following result tells us the degree of the resulting variety (or, to be precise, of the ideal generated by the minors) in case this bound is matched exactly. Here and everywhere else in the paper, by $\deg D_{k-1}(\matx M)$ we mean the degree of the \emph{ideal} generated by $k$-minors; if the minors define the determinantal variety in a reduced way, then this agrees with the degree of the variety.


\begin{theorem}[Thom--Porteous formula (corank one)~{\cite[Theorem 12.4]{EisenbudHarris2016}}]\label{thm:TP}
Let $0 < e \leq f$ be integers, and let $\matx{M}(\vect{x}) \in \Cset[\vect{x}]^{f\times e}$ be a matrix of homogeneous polynomials in $\vect{x} \in \Pset^n$, whose columns have constant degree $d_1,\dots,d_e$. 
Assume that each component $Z \subseteq D_{e-1}(\matx{M})$ of the determinantal variety of rank $e-1$ has the expected dimension $\dim Z = n - (f-e+1) \geq 0$. Then
\[
\deg D_{e-1}(\matx{M}) = T^{f-e+1}\text{-coefficient of }\frac{1}{(1-d_1T)\dotsm (1-d_eT)} \geq 0.
\]
\end{theorem}


If the dimension of $D_{e-1}(\matx{M})$ is zero, that is, if the determinantal variety consists of finitely many elements, then this gives a way to count the points $\vect{x} \in \Pset^n$ for which $\matx{M}(\vect{x})$ is of sub-maximal rank. 
This count, however, is only correct if we take the multiplicities into account, and is an overcount if the maximal minors of $\matx{M}$ do not define $D_{e-1}(\matx{M})$ in a reduced way.
A sufficient criterion for the maximal minors to define $D_{e-1}(\matx{M})$ in a reduced way is the following.

\begin{lemma}[{\cite[Proposition 0.4.15]{Leonie2025}}]\label{lem:reducedMinors}
Let $X \subseteq \Cset^n$ be an affine variety and consider a matrix polynomial $\matx{M}(\vect{x}) \in R_n^{f \times e}$, with $e \leq f$. 
Define a (quasi-projective variety $\mathcal{K}$ by the $f$ equations $\matx{M}(\vect{x})\vect{v} = \vect{0}$, that is,
\[
\mathcal{K} = \Set{ (\vect{x},[\vect{v}]) \in X \times \Pset^{e-1} | \matx{M}(\vect{x})\cdot \vect{v} = \vect{0}}.
\]
If $I(X)$ and these $f$ equations define $\mathcal K$ in a reduced way and if $\forall \, \vect{x} \in X: \rank \matx{M}(\vect{x}) \in \{e,e-1\}$, then the $e$-minors of $\matx{M}$ define $D_{e-1}(\matx{M})$ in a reduced way as well.
\end{lemma}

In this lemma, \enquote{defined in a reduced way}  is understood to be in the affine sense on each standard affine chart $\CC^{n+e-1} \cong \CC^n \times (\PP^{e-1}\setminus V(x_i)) \subseteq \CC^n\times \PP^{e-1}$.

\printbibliography[heading=bibintoc]

\end{document}